\documentclass[10pt]{article}
\usepackage{amsmath,amssymb,amsfonts,mathtools}
\usepackage{amsthm}
\usepackage{microtype}
\usepackage{geometry}
\usepackage{booktabs,tabularx,array,multirow}
\usepackage{enumitem}
\usepackage{tikz}
\usetikzlibrary{arrows.meta,positioning,calc,fit,shapes.geometric}
\usepackage{caption}
\usepackage[colorlinks=true,linkcolor=black,citecolor=black,urlcolor=black]{hyperref}
\setlist[enumerate]{itemsep=2pt,topsep=5pt}
\numberwithin{equation}{section}
\newtheorem{theorem}{Theorem}[section]
\newtheorem{proposition}[theorem]{Proposition}
\newtheorem{lemma}[theorem]{Lemma}
\newtheorem{corollary}[theorem]{Corollary}
\newtheorem{maintheorem}{Theorem}
\renewcommand{\themaintheorem}{\Alph{maintheorem}}
\theoremstyle{definition}
\newtheorem{definition}[theorem]{Definition}
\newtheorem{conjecture}[theorem]{Conjecture}
\theoremstyle{remark}
\newtheorem{remark}[theorem]{Remark}
\newcommand{\F}{\mathcal{F}}
\newcommand{\Ful}{\mathcal{F}}
\newcommand{\HH}{\mathcal{H}}
\newcommand{\FF}{\mathbb{F}}
\newcommand{\Z}{\mathbb{Z}}
\newcommand{\R}{\mathbb{R}}
\newcommand{\V}{V}
\newcommand{\E}{E}
\newcommand{\chd}{\chi_{d}}
\newcommand{\laml}{\lambda_{\ell}}
\newcommand{\scf}[1]{s_c(#1)}
\newcommand{\ind}{\mathbf{1}}
\DeclareMathOperator{\Alt}{Alt}
\DeclareMathOperator{\rank}{rank}

\title{From Diagonal Coloring to Proper Vertex Labelings\\Distinguishing Adjacent Facial Sums}
\author{Zhou Guofei\footnote{Email: gfzhou@nju.edu.cn}
\\
{\small School of Mathematics, Nanjing University, Nanjing 210093, China}
 }

\begin{document}
\maketitle
\begin{abstract}
In this paper we study proper vertex labelings of 3-connected plane graphs that give distinct boundary-vertex sums to faces sharing an edge. Let $\lambda(G)$ be the least $k$ for which labels from $\{0,\ldots,k\}$ suffice, and let $\laml(G)$ be its arbitrary-real-list counterpart. Our two main results are an eight-element list theorem for fullerenes and the universal bound $\lambda(G)\le708587$ for all 3-connected plane graphs. We establish structural nine- and eight-color bounds, classify triangulations whose vertices of degree at least five form an independent set, and determine their exact diagonal chromatic numbers. We also prove closure under triangular sums for $q\ge7$ diagonal colors, with an optimal uniform threshold. For cubic 3-connected plane graphs we obtain $\lambda\le14$, improving to $\lambda\le9$ in the bipartite case. 

\vskip 2mm
\noindent\textbf{Keywords.} diagonal coloring, facial sums, vertex labeling, plane triangulation, fullerene, list labeling\\
\noindent\textbf{2020 Mathematics Subject Classification.} 05C10, 05C15, 05C78.
\end{abstract}

\section{Introduction}\label{sec:intro}
Let $G$ be a plane graph and let $c$ assign numerical labels to its vertices. The \emph{facial sum} of a face $F$ is
\[
\scf{F}=\sum_{v\in V(F)}c(v).
\]
If two triangular faces $uvx$ and $uvy$ share an edge, then
\[
\scf{uvx}-\scf{uvy}=c(x)-c(y).
\]
Thus distinguishing the sums of edge-adjacent triangular faces is exactly the \emph{diagonal condition}: the vertices opposite their shared edge must receive different colors. For nontriangular faces, the corresponding difference remains linear but may involve several vertices.

For triangulations, cancellation converts the facial condition into a collection of pairwise coloring conflicts. For larger faces, the uncanceled labels give genuine linear inequalities. We therefore study the linear constraints directly while retaining properness on the original edges. This additional requirement matters: a labeling that distinguishes facial sums need not be proper, and a proper vertex coloring need not distinguish adjacent facial sums. Our goal is to satisfy both families of constraints using a small common integer alphabet or, where possible, small arbitrary real lists.

All graphs are finite and simple, and plane embeddings are viewed on the sphere. Every triangulation considered here has at least four vertices. In a 3-connected plane graph, facial boundaries are cycles, and two distinct faces meet in at most one vertex or one edge; see Diestel~\cite[Chapter~4]{diestel}. These intersection properties underpin the local and global counting arguments.

\subsection{Related work}
\begin{definition}\label{def:diagonal}
In a plane triangulation $T$, vertices $x,y$ are \emph{diagonally adjacent} if $uvx$ and $uvy$ are the two faces incident with some edge $uv$. A \emph{diagonal coloring} is a proper vertex coloring that separates every diagonally adjacent pair. Its minimum number of colors is $\chd(T)$.
\end{definition}
Bouchet, Fouquet, Jolivet, and Rivi\`ere~\cite{bouchet} introduced diagonal coloring, proved a twelve-color bound, and constructed a triangulation requiring nine colors. Borodin~\cite{borodin1990} improved the upper bound to eleven, and Sanders and Zhao~\cite{sanders1995} reduced it to ten. Thus the classical bounds are
\[
9\le \sup_T\chd(T)\le10.
\]
The nine-color conjecture asks whether the upper bound is nine. Huang et al.~\cite{huang} proved this for triangulations of diameter at most two and characterized the equality case. The stellated prism used below is the corresponding nine-vertex example.

These results concern the ordinary diagonal chromatic number. The degenerate diagonal colorings of Kierstead, Mohar, \v{S}pacapan, Yang, and Zhu~\cite{kierstead} impose additional restrictions, so their eleven-color and twelve-list-color bounds are not competing estimates for the same minimum. More generally, $d$-diagonal coloring separates vertices that can become cofacial after deletion of at most $d$ edges~\cite{sanders1996,borodin1999}. Kr\'al' and \v{S}krekovski~\cite{kral} disproved a sharpness conjecture for diagonal coloring on surfaces of positive Euler genus; that conjecture is distinct from the plane nine-color conjecture.

Facial numerical weightings have a separate history. Face-antimagic labelings and arithmetic-progression variants were studied by Ba\v{c}a~\cite{baca2005}; the survey of Ba\v{c}a et al.~\cite{bacasurvey} covers their vertex-, edge-, and face-labeling versions. Face-irregular evaluations distinguish \emph{all} face weights; see Ba\v{c}a et al.~\cite{bacairregular} and Bayati et al.~\cite{bayati} for grid families. The labeled elements, permitted alphabets, and required face comparisons differ from those used here. In particular, these conditions need not include properness of vertex labels. Entire face-coloring weightings allow labels on vertices, edges, and faces~\cite{jendrol}. Our constraints \eqref{eq:proper} and \eqref{eq:facial} instead use only vertex labels and require both properness and local facial-sum distinction.

\subsection{Main results and organization}
We seek labelings that are proper on vertices and distinguish the sums of edge-adjacent faces. Write $\lambda(G)$ for the least $k$ such that the common alphabet $\{0,\ldots,k\}$ suffices, and $\laml(G)$ for the least $k$ such that every assignment of $(k+1)$-element real lists admits such a labeling. Thus a bound of $k$ means that $k+1$ labels are available. Section~\ref{sec:prelim} gives the formal definitions.

Our principal existence results address two different settings. A \emph{fullerene} is a 3-connected cubic plane graph whose faces have length five or six.

\begin{maintheorem}[Fullerene lists; Theorem~\ref{thm:fullerene-list}]\label{thm:intro-fullerene}
Every fullerene admits a proper vertex labeling that distinguishes edge-adjacent facial sums from arbitrary eight-element real lists. In particular,
\[
\lambda(G)\le\laml(G)\le7.
\]
\end{maintheorem}
The proof first labels the vertices outside a pentagonal face greedily, then extends the labeling across that face using a polynomial in five variables. A single coefficient controls every combination of neighboring pentagons and hexagons. Appendix~\ref{app:coefficient} proves its nonvanishing symbolically over $\FF_7$.

\begin{maintheorem}[Universal bound; Proposition~\ref{prop:edge-count} and Theorem~\ref{thm:universal}]\label{thm:intro-universal}
Every 3-connected plane graph $G$ satisfies
\[
\lambda(G)\le\min\{3|V(G)|-6,\;708587\}.
\]
\end{maintheorem}
The constant $708587=4\cdot3^{11}-1$ is independent of the order, degrees, and face lengths of $G$. Its proof bounds the number of projected facial-difference types linearly in the number of selected vertices, applies matroid partition, and constructs a labeling by successive modular digits. This argument gives a common integer alphabet; it does not establish a constant bound for arbitrary lists. For such lists, Proposition~\ref{prop:edge-count} gives the separate estimate $\laml(G)\le|E(G)|$.

The triangulation results form a complementary structural part of the paper. Section~\ref{sec:structural} establishes nine diagonal colors when all-high facial triangles are edge-disjoint and eight colors when no such triangle exists; here \emph{high} means degree at least five in the original triangulation. We also classify triangulations with independent high vertices, determine their exact diagonal chromatic numbers, and prove that triangular sums preserve diagonal $q$-colorability for $q\ge7$, with an optimal uniform threshold. These results concern structural subclasses and do not resolve the unrestricted nine-color conjecture. Their proofs provide the diagonal-coloring viewpoint and examples used to place the later labeling bounds in context; the proofs of Theorems~\ref{thm:intro-fullerene} and~\ref{thm:intro-universal} use separate constructions.

Between the triangulation and fullerene cases, Section~\ref{sec:cubic} treats cubic 3-connected plane graphs through a dual incidence map. It proves $\lambda\le14$, improves this to $\lambda\le9$ for bipartite graphs, and gives the twelve-label bound $\lambda\le11$ in the modular cases specified in Theorem~\ref{thm:cubic}.

After the preliminaries in Section~\ref{sec:prelim}, the paper proceeds from structural diagonal coloring in Section~\ref{sec:structural} to cubic graphs in Section~\ref{sec:cubic}, fullerenes in Section~\ref{sec:fullerenes}, and the universal construction in Section~\ref{sec:universal}. Section~\ref{sec:outlook} collects the bounds, open questions, and verification scope. Appendix~\ref{app:coefficient} contains the coefficient proof, while Appendix~\ref{app:elementary} gives a weaker elementary fullerene bound. The Four Color Theorem is used in specified results; all imported theorem-level ingredients are identified where they enter.

\section{Labeling parameters and preliminary tools}\label{sec:prelim}
\subsection{Parameters and conventions}
\begin{definition}\label{def:label}
An \emph{adjacent-facial-sum-distinguishing proper vertex labeling} of a 3-connected plane graph $G$ is a map $c:V(G)\to\R$ satisfying
\begin{align}
c(u)&\ne c(v) &&(uv\in E(G)),\label{eq:proper}\\
\scf{F}&\ne\scf{H} &&(F,H\text{ share an edge}).\label{eq:facial}
\end{align}
Let $\lambda(G)$ be the least integer $k\ge0$ for which such a labeling takes values in $\{0,1,\ldots,k\}$. Let $\laml(G)$ be the least $k\ge0$ such that every assignment of real lists $A(v)$ with $|A(v)|=k+1$ admits a labeling satisfying \eqref{eq:proper}--\eqref{eq:facial} and $c(v)\in A(v)$.
\end{definition}
Plainly $\lambda(G)\le\laml(G)$. We also write
\[
K_3(n)=\max_{\substack{|V(G)|=n\\G\ \text{3-connected plane}}}\lambda(G),\quad
K_3^*=\sup_{n\ge4}K_3(n),\quad
K_{\Ful}=\sup_{G\in\Ful}\lambda(G),
\]
where $\Ful$ denotes the family of fullerenes.

\begin{remark}\label{rem:alphabet}
The zero-based alphabet is part of the definition. Unlike proper coloring, facial-sum distinction is not generally invariant under translation. For the translated labeling $(c+t)(v)=c(v)+t$,
\[
s_{c+t}(F)-s_{c+t}(H)=\scf{F}-\scf{H}+t(|F|-|H|).
\]
Translation preserves the requirement when the compared faces have the same length. We permute final vertex color names only for triangulations, where the constraints are pairwise inequalities. Later permutations of \emph{auxiliary} dual colors occur before the primal face-sum congruences are solved; no invariance is claimed for a completed primal labeling.
\end{remark}

\subsection{The diagonal-coloring correspondence}
\begin{proposition}\label{prop:correspondence}
For every plane triangulation $T$,
\[
\lambda(T)=\chd(T)-1.
\]
Moreover, $\laml(T)$ is the list chromatic number of the diagonal conflict graph minus one.
\end{proposition}
\begin{proof}
For two triangular faces $uvx$ and $uvy$, the difference of facial sums is $c(x)-c(y)$. Thus the edge constraints and facial comparisons are exactly the proper-coloring constraints of the graph obtained by adding all diagonal conflicts to $T$. For ordinary coloring, colors may be renamed by consecutive integers. For arbitrary real lists, the same equivalence is precisely list coloring of that conflict graph.
\end{proof}
In particular, if $K_{\mathrm{tri}}=\sup_T\lambda(T)$, the classical diagonal bounds give
\begin{equation}\label{eq:tri-bounds}
8\le K_{\mathrm{tri}}\le9,\qquad K_3^*\ge8.
\end{equation}

\subsection{Polynomial bounds and a modular combination lemma}
For each edge $e=uv$, choose an order of its endpoints and an order $F_e,H_e$ of its incident faces. Define
\begin{equation}\label{eq:PG}
P_G(x)=\prod_{e=uv\in E(G)}(x_u-x_v)\bigl(S_{F_e}(x)-S_{H_e}(x)\bigr),\qquad
S_F(x)=\sum_{v\in V(F)}x_v.
\end{equation}
The boundary vertex sets of adjacent faces differ, so each linear factor is nonzero. Consequently $P_G$ is a nonzero polynomial over $\R$, and its nonzero evaluations are exactly the desired labelings. We use the standard \emph{grid principle}: if a nonzero polynomial has degree at most $D_v$ in each variable $x_v$, it does not vanish on a product of real sets of cardinalities $D_v+1$. This follows by successive univariate interpolation and is a basic polynomial-method consequence of the Combinatorial Nullstellensatz~\cite{alon}.

\begin{proposition}[Edge-count list bound]\label{prop:edge-count}
If $G$ is a 3-connected plane graph with $n$ vertices and $m$ edges, then
\[
\lambda(G)\le\laml(G)\le m\le3n-6.
\]
\end{proposition}
\begin{proof}
Fix $v$ and an edge $e=uw$. If $v$ is an endpoint of $e$, its variable occurs in $x_u-x_w$ but cancels from $S_{F_e}-S_{H_e}$, since both incident faces contain $v$. If $v$ is not an endpoint, it can occur only in the facial factor. Thus the pair of factors indexed by an edge has degree at most one in $x_v$, and $\deg_{x_v}P_G\le m$. The grid principle gives a nonzero evaluation from arbitrary $(m+1)$-element real lists. The planar edge bound completes the proof.
\end{proof}

\begin{proposition}[Local list bound]\label{prop:local-list}
For $v\in V(G)$, put
\[
D_v=d_G(v)+\sum_{F\ni v}(|F|-2).
\]
For any assignment of lists satisfying $|A(v)|\ge\min\{m,D_v\}+1$ at every vertex, there is a labeling as in Definition~\ref{def:label}. In particular, if $\Delta(G)\le\Delta$ and all faces have length at most $r$, then
\[
\laml(G)\le\min\{m,\Delta(r-1)\}.
\]
\end{proposition}
\begin{proof}
For a face $F$ containing $v$, comparisons across the two boundary edges of $F$ incident with $v$ cancel $x_v$. The remaining $|F|-2$ edges account for every facial factor in which $x_v$ occurs, possibly with overcounting. Hence $\deg_{x_v}P_G\le D_v$. Proposition~\ref{prop:edge-count} also gives $\deg_{x_v}P_G\le m$. Apply the grid principle with the smaller bound, then use $D_v\le d_G(v)(r-1)$.
\end{proof}

\begin{lemma}[Modular combination]\label{lem:modular-combination}
Suppose $G$ has a proper $r$-coloring and a map $b:V(G)\to\{0,\ldots,q-1\}$ for which the sums of edge-adjacent faces differ modulo $q$. Then $\lambda(G)\le rq-1$.
\end{lemma}
\begin{proof}
Take a proper $a:V(G)\to\{0,\ldots,r-1\}$ and set $c=qa+b$. On an edge $uv$, the nonzero multiple $q(a(u)-a(v))$ cannot be canceled by $b(u)-b(v)$, whose absolute value is less than $q$. Facial differences of $c$ agree modulo $q$ with those of $b$. Finally, $0\le c(v)\le rq-1$.
\end{proof}

\section{Structural diagonal-coloring bounds}\label{sec:structural}
Let $D(T)$ be the graph on $V(T)$ whose edges are the original edges of $T$ together with all diagonal conflicts. Thus
\begin{equation}\label{eq:conflict}
\chi(D(T))=\chd(T)=\lambda(T)+1.
\end{equation}
Set
\[
H=\{v\in V(T):d_T(v)\ge5\},\qquad L=V(T)\setminus H,\qquad \HH=T[H],
\]
and let $\F_H$ denote the set of facial triangles whose three vertices lie in $H$. The term \emph{edge-disjoint all-high facial triangles} means that no two members of $\F_H$ share an edge; they may share vertices. The qualifier ``facial'' is essential: a separating triangle on high vertices need not belong to $\F_H$.

\begin{lemma}\label{lem:independent-planar}
If $I$ is independent in $T$, then $D(T)[I]$ is planar.
\end{lemma}
\begin{proof}
For a diagonal conflict $xy$ witnessed by the faces $uvx,uvy$, draw the edge $xy$ through the midpoint of $uv$, within the interiors of those two faces. Since $I$ is independent, no facial triangle contains two vertices of $I$; hence at most one such drawn segment is needed within any face. The segments have no crossings. There are no original edges inside $I$, and duplicate diagonal edges may be deleted.
\end{proof}

\begin{lemma}\label{lem:degree-conflict}
Every vertex satisfies $d_{D(T)}(v)\le2d_T(v)$. In particular, vertices of $L$ have conflict degree at most eight.
\end{lemma}
\begin{proof}
In addition to its original neighbors, $v$ has at most one opposite vertex across each edge of its link cycle. Repetitions only decrease the number of distinct conflicts.
\end{proof}

\subsection{Nine and eight colors}
The following construction converts conflicts among high vertices into the edges of a $1$-planar supergraph. It is the geometric ingredient in both color bounds.

\begin{figure}[t]
\centering
\begin{tikzpicture}[scale=0.9, font=\footnotesize]
  \foreach \shift/\name in {0/{(a) Low-vertex star},5.0/{(b) Clique replacement}} {
    \begin{scope}[xshift=\shift cm]
      \coordinate (a) at (0,1.4); \coordinate (b) at (1.4,0); \coordinate (c) at (0,-1.4); \coordinate (d) at (-1.4,0);
      \ifdim\shift pt<1pt
        \draw[gray!55,dashed,thin] (a)--(b)--(c)--(d)--cycle;
      \else
        \draw[thick] (a)--(b)--(c)--(d)--cycle;
      \fi
      \ifdim\shift pt<1pt
        \draw[thick] (0,0)--(a) (0,0)--(b) (0,0)--(c) (0,0)--(d);
        \filldraw[fill=white,draw=black] (0,0) circle (0.14) node[below right] {$v$};
      \else
        \draw[thick] (a)--(c) (b)--(d);
      \fi
      \foreach \pt/\lab/\pos in {a/{$h_1$}/above,b/{$h_2$}/right,c/{$h_3$}/below,d/{$h_4$}/left} {
         \filldraw[fill=black,draw=black] (\pt) circle (0.09);
         \node[\pos] at (\pt) {\lab};
      }
      \node at (0,-2.2) {\name};
    \end{scope}
  }
\end{tikzpicture}
\caption{Replacing the star of a degree-four low vertex by a clique on its high neighbors. The dashed cycle in (a) only indicates terminal order, not original graph edges. In (b) only the two diagonals cross, and each resulting edge is crossed at most once. Degree-two and degree-three replacements are planar.}\label{fig:oneplanar}
\end{figure}
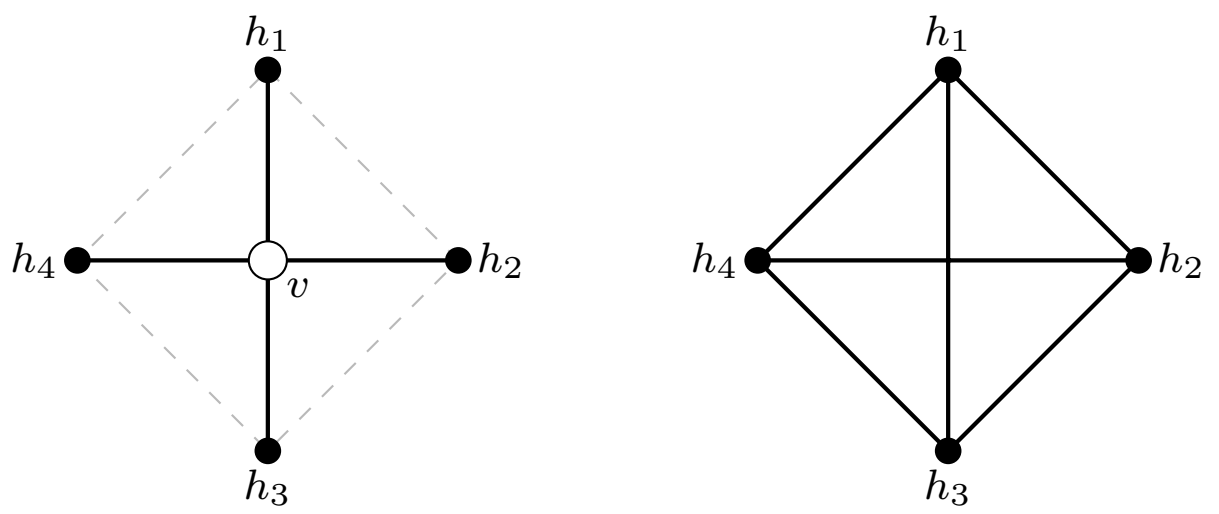

\begin{lemma}\label{lem:high-oneplanar}
If no two faces in $\F_H$ share an edge, then $D(T)[H]$ is a subgraph of a $1$-planar graph.
\end{lemma}
\begin{proof}
Delete each edge with both endpoints in $L$ and subdivide every edge with both endpoints in $H$. The resulting plane graph is bipartite: one part is $H$, and the other consists of original low vertices and subdivision vertices. Every vertex of the latter part has degree at most four.

Replace each star centered in the latter part by a clique on its high neighbors, routing edges along parallel strands in thin strips around the original incidences and inside a small disk about the center. The disks and strips can be chosen with pairwise disjoint interiors. A clique on at most three terminals is planar. For four terminals in cyclic order, a $K_4$ can be drawn with one crossing between the two diagonals, as in Figure~\ref{fig:oneplanar}. No new edge is crossed more than once. After deletion of duplicate edges, denote the resulting $1$-planar graph by $P$.

Every original high--high edge is represented in $P$ through its subdivision vertex. Suppose $xy$ is a diagonal conflict with $x,y\in H$, witnessed by $uvx$ and $uvy$. If both $u$ and $v$ were high, those two facial triangles would share the all-high edge $uv$, contrary to the hypothesis. Thus an endpoint of $uv$ is low and has both $x$ and $y$ as high neighbors. Its star replacement therefore supplies $xy$. Hence $D(T)[H]\subseteq P$.
\end{proof}

\begin{theorem}\label{thm:nine}
If no two faces in $\F_H$ share an edge, then
\[
\chd(T)\le9,\qquad\lambda(T)\le8.
\]
\end{theorem}
\begin{proof}
Lemma~\ref{lem:high-oneplanar} and Borodin's six-color theorem for $1$-planar graphs~\cite{borodin1984,borodin1995} give a six-coloring of $D(T)[H]$. Extend it greedily to $L$ using a palette of nine colors: each vertex of $L$ has at most eight conflict neighbors by Lemma~\ref{lem:degree-conflict}.
\end{proof}

\begin{theorem}\label{thm:eight}
If $\F_H=\varnothing$, then
\[
\chd(T)\le8,\qquad\lambda(T)\le7.
\]
\end{theorem}
\begin{proof}
\emph{Step 1: residual lists and components with slack.}
By Lemma~\ref{lem:high-oneplanar} and Borodin's six-color theorem~\cite{borodin1984,borodin1995}, color $D(T)[H]$ from $\{0,\ldots,5\}$. For $v\in L$, let $A(v)$ consist of the colors in $\{0,\ldots,7\}$ unused by its high conflict neighbors, and set $Q=D(T)[L]$. Then
\begin{equation}\label{eq:residual}
|A(v)|\ge8-d_{D(T)}(v,H)\ge d_Q(v).
\end{equation}
Consider each component $C$ of $Q$. If some $r\in V(C)$ has total conflict degree at most seven, then $|A(r)|\ge d_C(r)+1$. Root a spanning tree of $C$ at $r$ and process vertices in reverse tree order. Each nonroot vertex still has its uncolored parent when it is processed, while the root has one extra available color. This gives a list coloring of $C$.

\emph{Step 2: components with no slack.}
Otherwise every vertex $v$ of $C$ has total conflict degree eight. By Lemma~\ref{lem:degree-conflict}, $d_T(v)=4$. Write its four original neighbors in cyclic link order as $a_0,\ldots,a_3$, and let $w_i$ be opposite $v$ across the edge $a_i a_{i+1}$ (indices modulo four). Equality of the total conflict degree to eight forces $w_0,\ldots,w_3$ to be distinct and outside $\{v,a_0,\ldots,a_3\}$.

Each $a_i$ has the five distinct neighbors $v,a_{i-1},a_{i+1},w_{i-1},w_i$, so $a_i\in H$. Since $a_i a_{i+1}w_i$ is facial and no facial triangle is all-high, every $w_i$ belongs to $L$. Thus every vertex of $C$ has its four original neighbors in $H$ and its four diagonal neighbors in $L$, necessarily in the same component $C$. The graph $C$ is therefore four-regular and its vertices are pairwise nonadjacent in $T$. By Lemma~\ref{lem:independent-planar}, $C$ is planar, hence cannot be $K_5$.

The list form of Brooks' theorem~\cite{erdosrubin} states that a connected graph is colorable from lists of size equal to its maximum degree unless it is complete or an odd cycle. Neither exception applies to the four-regular planar graph $C$. Its residual lists have size at least four by \eqref{eq:residual}, so $C$ can be colored. Combining the components completes the extension.
\end{proof}

\begin{corollary}\label{cor:bipartite-high}
If $\HH$ is bipartite, then $\chd(T)\le8$. Equivalently, the conclusion holds if every odd cycle of $T$ contains a vertex of degree three or four.
\end{corollary}
\begin{proof}
A bipartite $\HH$ contains no facial triangle. The equivalence follows from the odd-cycle characterization of bipartiteness.
\end{proof}

\subsection{Independent high vertices: exact classification}
Let $B_m$ be the cycle bipyramid obtained from a cycle $C_m$ by joining two \emph{nonadjacent} apices to every cycle vertex, where $m\ge3$.

\begin{theorem}[Structural classification]\label{thm:classification}
The set $H$ is independent in $T$ if and only if $T\cong K_4$ or $T\cong B_m$ for some $m\ge3$.
\end{theorem}
\begin{proof}
Write $h=|H|$ and $\ell=|L|$. As every face contains at most one high vertex,
\[
\sum_{v\in H}d_T(v)\le 2(h+\ell)-4.
\]
Together with $d_T(v)\le4$ for $v\in L$ and Euler's degree identity, this gives
\[
6(h+\ell)-12\le 2(h+\ell)-4+4\ell,
\]
so $h\le2$.

If $h=2$, equality holds throughout: every low vertex has degree four, and every facial triangle contains exactly one high vertex. Let $v\in L$ and list its four distinct neighbors in cyclic order. Since consecutive neighbors form a face with $v$, exactly one vertex from each consecutive pair is high. The link therefore alternates between high and low vertices. Thus $v$ has the two vertices of $H$ as its high neighbors and exactly two low neighbors. Every vertex of $L$ is adjacent to both high vertices. The link of either high vertex is a single cycle containing all of $L$; moreover, each low vertex already has its two low neighbors on this cycle. Hence $T\cong B_\ell$.

If $h=1$, let $x$ be the high vertex and set $n=|V(T)|$. Using $d_T(x)\le n-1$ and $d_T(v)\le4$ otherwise in the degree sum gives $n\le7$. Since $d_T(x)\ge5$, also $n\ge6$. For $n=6$ the degree sum forces $d_T(x)=5$, and for $n=7$ it forces $d_T(x)=6$; thus $x$ is universal. Deleting it leaves a maximal outerplane graph on $n-1$ vertices, which has at least two degree-two ears (the triangle-adjacency tree has at least two leaves). These ears have degree three in $T$. Consequently
\[
6n-12\le(n-1)+4(n-1)-2=5n-7,
\]
which implies $n\le5$, a contradiction.

If $h=0$, the degree sum gives $6n-12\le4n$, whence $4\le n\le6$. These cases yield $K_4$, $B_3$, and $B_4$, respectively: on five vertices the missing edge identifies the two apices, and on six vertices all degrees are four, so the complement is a perfect matching and $T$ is the octahedron.

Conversely, $K_4$ has no high vertices, and the only possible high vertices of $B_m$ are its two nonadjacent apices. Thus all the listed graphs have independent high vertices.
\end{proof}

\begin{proposition}[Exact diagonal chromatic numbers]\label{prop:bipyramid-colors}
For $m\ge3$,
\begin{equation}\label{eq:bipyramid}
\chd(K_4)=4,\qquad
\chd(B_m)=
\begin{cases}
5,&3\mid m,\\
7,&m=5,\\
6,&3\nmid m\text{ and }m\ne5.
\end{cases}
\end{equation}
In particular, triangulations with independent high vertices satisfy $\chd\le7$, equivalently $\lambda\le6$, with equality precisely for $B_5$.
\end{proposition}
\begin{proof}
Both apices conflict with every cycle vertex and with each other. Two cycle vertices conflict exactly when their cyclic distance is at most two: original cycle edges provide distance-one conflicts, and faces across an apex--cycle edge provide distance-two conflicts. Therefore
\[
D(B_m)=K_2\vee C_m^2,\qquad \chd(B_m)=2+\chi(C_m^2),
\]
where $\vee$ denotes the graph join and $C_m^2$ is the square of the cycle.

For $m=3,4,5$, the squares are $K_3,K_4,K_5$. For $m\ge6$, the square is connected, four-regular, and not $K_5$; Brooks' theorem~\cite{brooks} gives $\chi(C_m^2)\le4$. In a proper three-coloring of $C_m^2$, each three consecutive cycle vertices use all three colors, so the coloring repeats with period three. It exists exactly when $3\mid m$, and repeating three colors is sufficient in those cases. Four colors are otherwise necessary, except when $C_5^2=K_5$. Finally, $D(K_4)=K_4$.
\end{proof}

\subsection{Independent odd-cycle transversals}
\begin{theorem}\label{thm:transversal}
If $S\subseteq H$ is independent in $D(T)$ and $T[H\setminus S]$ is bipartite, then $\chd(T)\le9$.
\end{theorem}
\begin{proof}
Let $A\cup B$ be the bipartition of $T[H\setminus S]$. The sets $A$ and $B$ are independent in $T$, so Lemma~\ref{lem:independent-planar} makes both $D(T)[A]$ and $D(T)[B]$ planar. Color them with disjoint palettes of four colors using the Four Color Theorem~\cite{fourcolor}. Give the independent set $S$ a ninth color and extend greedily over $L$, using Lemma~\ref{lem:degree-conflict}.
\end{proof}

\begin{corollary}\label{cor:obstructions}
If $\chd(T)\ge10$, then two faces in $\F_H$ share an edge and $T[H\setminus\{v\}]$ is nonbipartite for every $v\in H$.
\end{corollary}
\begin{proof}
Apply Theorems~\ref{thm:nine} and~\ref{thm:transversal} contrapositively, with $S=\{v\}$ in the latter.
\end{proof}

\subsection{Triangular sums and the optimal threshold}
A \emph{triangular sum} identifies a facial triangle of each of two plane triangulations and removes the two identified face interiors. New diagonal conflicts may arise across the three edges of the common triangle.

\begin{proposition}[Closure under triangular sums]\label{prop:gluing}
For every $q\ge7$, a plane triangular sum of two diagonally $q$-colorable triangulations is diagonally $q$-colorable.
\end{proposition}
\begin{proof}
Let $abc$ be the identified facial triangle. Permute the color names on one piece to match its three boundary colors with those on the other. On either side of the seam, the vertex opposite a boundary edge uses a color outside the three boundary colors: properness separates it from the endpoints, and the original diagonal constraint separates it from the third boundary vertex.

Only the three seam edges can create new diagonal conflicts. The same pair of vertices may occur across more than one seam edge, so there are at most three distinct new conflicts. Randomly and uniformly permute the remaining $q-3$ colors on one piece, fixing the three boundary colors. For each seam edge, the probability that its two opposite vertices receive the same color is $1/(q-3)$. By the union bound, the probability of at least one failure is at most $3/(q-3)<1$ when $q\ge7$. A permutation avoiding all failures therefore exists.
\end{proof}

\begin{proposition}[Optimality of the uniform gluing threshold]\label{prop:gluing-sharp}
The assumption $q\ge7$ in Proposition~\ref{prop:gluing} cannot be replaced by $q\ge6$. A triangular sum of $K_4$ and the octahedron has diagonal conflict graph $K_7$, although its two summands have diagonal chromatic numbers four and six.
\end{proposition}
\begin{proof}
Let $abc$ be the identified face and $x$ the fourth vertex of $K_4$. The octahedron has diagonal conflict graph $K_6$. Its three vertices outside $abc$ are exactly those opposite the seam edges $ab,bc,ca$ on the octahedral side. In the sum, $x$ is adjacent to $a,b,c$ and is diagonally opposite the other three vertices across the seam edges.

All conflicts among the six octahedral vertices survive. Indeed, the octahedron is $K_{2,2,2}$, and each of its three nonadjacent pairs is diagonally opposite across all four edges joining the other two parts. Removing the face $abc$ destroys only one of these four witnesses for each pair. The six original vertices therefore still induce $K_6$, and $x$ conflicts with all six. Hence the full conflict graph is $K_7$.
\end{proof}

\begin{corollary}\label{cor:separating}
For $q\ge7$, if every piece obtained by recursively cutting along separating triangles is diagonally $q$-colorable, then so is $T$. In particular, if a triangulation is not diagonally nine-colorable, at least one of its pieces without separating triangles also fails to be nine-colorable.
\end{corollary}
\begin{proof}
Reverse the cuts and apply Proposition~\ref{prop:gluing} successively.
\end{proof}

\subsection{Sharpness and examples}
\begin{proposition}\label{prop:stellated-prism}
Insert a degree-four center into each quadrilateral face of a triangular prism and join it to the four boundary vertices. The resulting triangulation $T_{\mathrm{pr}}$ satisfies $D(T_{\mathrm{pr}})\cong K_9$. Its only all-high facial triangles are the two prism ends. Thus Theorem~\ref{thm:nine} is sharp.
\end{proposition}
\begin{proof}
The six original vertices have degree five, and the three inserted vertices have degree four. Every nonadjacent pair of original vertices is opposite across an edge incident with an inserted center. An inserted center and an original vertex outside its quadrilateral are opposite across an edge of a prism end. Finally, two inserted centers are opposite across the common prism edge of their quadrilaterals. These cases account for all nonedges, so the conflict graph is $K_9$.

The prism ends remain the only all-high facial triangles and are edge-disjoint. Thus this example meets the hypothesis of Theorem~\ref{thm:nine}, but not the stronger hypothesis of Theorem~\ref{thm:eight}.
\end{proof}

The barycentric subdivision of any 3-connected plane graph satisfies Theorem~\ref{thm:eight}: each facial triangle contains an edge-vertex of degree four. The same is true after inserting a degree-four center into every face of a plane quadrangulation. Conversely, the icosahedron has no low vertices and satisfies neither structural hypothesis, although its conflict graph is $K_{12}$ minus a perfect matching and is six-colorable.

\section{Cubic 3-connected plane graphs}\label{sec:cubic}
Let $G$ be cubic and 3-connected, and let $J=G^*$ be its plane dual. Write the face--vertex incidence matrix as
\[
B\in\{0,1\}^{\F(G)\times V(G)},\qquad B_{F,v}=\ind_{\{v\in V(F)\}},
\]
where $\F(G)$ denotes the set of faces. A primal vertex $v$ lies on three faces, whose corresponding dual vertices form a facial triangle of $J$. Thus the column $B_{\cdot,v}$ is the indicator of this triangle in $\Z^{V(J)}=\Z^{\F(G)}$. In particular, $Bc$ is precisely the vector of primal facial sums. We first describe the integral obstruction to prescribing that vector, then choose suitable proper dual colors that satisfy the obstruction.

\subsection{The integral obstruction to prescribed face sums}
Declare two vertices of $J$ equivalent if they are opposite vertices of adjacent facial triangles, and take the transitive closure. Denote the equivalence classes by $C_1,\ldots,C_s$. Two columns corresponding to adjacent dual triangles differ by $e_c-e_d$ for equivalent vertices $c,d$. Hence all facial triangles have the same vector of class counts. Every dual vertex belongs to a triangle, so $1\le s\le3$, and the count vectors, after reordering, are $(3)$, $(1,2)$, and $(1,1,1)$.

\begin{lemma}\label{lem:cokernel}
The cokernel $\Z^{\F(G)}/B\Z^{V(G)}$ is, respectively, $\Z/3\Z$, $\Z$, or $\Z^2$ in the three cases above. For any integer $q\ge1$, the congruence $Bb=h\pmod q$ is solvable precisely under the corresponding condition:
\begin{enumerate}[label=\textup{(\arabic*)}]
\item $s=1$: $\sum_{F}h(F)\in3(\Z/q\Z)$;
\item $s=2$: $\sum_{F\in D}h(F)=2\sum_{F\in A}h(F)\pmod q$, where classes $A,D$ occur once and twice, respectively, in each triangle;
\item $s=3$: the sums of $h$ over the three equivalence classes are equal modulo $q$.
\end{enumerate}
\end{lemma}
\begin{proof}
Differences of adjacent triangle columns generate all differences of standard basis vectors within an equivalence class. Quotienting $\Z^{V(J)}$ by these differences identifies all basis vectors within each class, leaving $\Z^s$ with the class sums as coordinates. Every triangle column then has the same image $a$, where $a=(3)$, $(1,2)$, or $(1,1,1)$. Consequently the original cokernel is $\Z^s/\langle a\rangle$, yielding the three groups in the statement. The same presentation reduced modulo $q$ gives the stated congruence criteria.
\end{proof}

\subsection{Balancing independent color classes}
All weights in this subsection are residues, while cardinalities of vertex classes are ordinary integers. We need a proper coloring of the dual satisfying one linear congruence.

\begin{lemma}[Four nonzero class weights]\label{lem:four-weights}
If $W_0,W_1,W_2,W_3\in\FF_5^\times$ sum to zero and no pair sums to zero, then their multiset is one of
\[
(1,1,1,2),\quad (1,3,3,3),\quad (2,2,2,4),\quad (3,4,4,4).
\]
\end{lemma}
\begin{proof}
The pairs of opposite nonzero residues in $\FF_5$ are $\{1,4\}$ and $\{2,3\}$. Thus the support of the multiset contains at most one member of each opposite pair. A singleton support $\{a\}$ would give $4a\ne0$, so the support is exactly one of $\{1,2\}$, $\{1,3\}$, $\{2,4\}$, or $\{3,4\}$. If the first value is $a$ and occurs $j$ times, solve $aj+b(4-j)=0$ in $\FF_5$ with $1\le j\le3$. The respective values $j=3,1,3,1$ give the four listed multisets.
\end{proof}

\begin{lemma}[Balancing modulo five]\label{lem:balance-five}
Let $J$ be properly four-colorable with $|V(J)|\ge5$, and let $w(v)\in\{1,3\}\subseteq\FF_5$. There is a proper map $h:V(J)\to\FF_5$ such that
\[
\sum_{v\in V(J)}w(v)h(v)=0.
\]
\end{lemma}
\begin{proof}
Fix a proper four-coloring with independent color classes $V_0,\ldots,V_3$, allowing empty classes. We may split one class into independent subclasses and give distinct colors to all nonempty resulting classes; at most five classes are used. Let
\[
W_i=\sum_{v\in V_i}w(v),\qquad W=\sum_{i=0}^3W_i,\qquad L=\sum_{i=0}^3iW_i.
\]

\emph{Step 1: translation and an unused color.}
If $W\ne0$, assigning the color $i+t$ to $V_i$ changes the weighted sum to $L+tW$. Choose $t=-L/W$, so assume henceforth
\begin{equation}\label{eq:Wzero}
W=0.
\end{equation}
Initially use colors $0,1,2,3$, leaving color $4$ unused. If a class $V_i$ contains at least four vertices, choose four of them. Any subset may be recolored to $4$ without losing properness. Each chosen vertex contributes a nonzero increment $d_v=(4-i)w(v)$ to the weighted sum.

The subset sums of four nonzero elements of $\FF_5$ fill $\FF_5$: starting from $U=\{0\}$, adjoining the option to use an increment $d\ne0$ replaces $U$ by $U\cup(U+d)$. If $U$ is a nonempty proper subset of $\FF_5$, then $U+d\ne U$, because $d$ generates its additive group. Thus the subset-sum set grows by at least one at each step until it is the whole field. Four increments suffice, so we can realize $-L$ and obtain weighted sum zero.

We may therefore assume that $|V_i|\le3$ for every $i$. The possible multisets of vertex weights in a nonempty class are
\begin{center}
\begin{tabular}{c l}
\toprule
Class size & Vertex-weight multisets \\
\midrule
1 & $(1),\ (3)$\\
2 & $(1,1),\ (1,3),\ (3,3)$\\
3 & $(1,1,1),\ (1,1,3),\ (1,3,3),\ (3,3,3)$\\
\bottomrule
\end{tabular}
\end{center}
In particular, total weight two comes from $(1,1)$ or $(1,3,3)$; weight four from $(1,3)$ or $(3,3,3)$; weight one from a singleton or $(3,3)$; and weight three from a singleton or $(1,1,1)$.

\emph{Step 2: at most three nonzero class weights.}
If every $W_i$ is zero, any assignment of distinct colors works. Exactly one nonzero $W_i$ contradicts \eqref{eq:Wzero}. If exactly three are nonzero, call them $a,b,c$, with $a+b+c=0$. Give their classes colors $0,1,-b/c$. These are distinct: $-b/c$ is neither $0$ nor $1$, since $b,c,a$ are nonzero. Zero-weight classes receive distinct unused colors. This works for at most five total classes, provided exactly three have nonzero weights.

If exactly two weights are nonzero, they are $a,-a$, namely either $\{1,4\}$ or $\{2,3\}$. In the first case split a vertex off the weight-four class, whose pattern is $(1,3)$ or $(3,3,3)$; in the second case split a vertex off the weight-two class, whose pattern is $(1,1)$ or $(1,3,3)$. Both resulting subclasses have nonzero weights. We now have exactly three nonzero weights among at most five independent classes, so the preceding construction applies.

\emph{Step 3: four nonzero class weights.}
If an opposite pair occurs, the weights can be ordered as $(a,-a,b,-b)$. Assign the first two classes colors $0,1$. The oriented differences of distinct members of $\{2,3,4\}$ cover $\FF_5^\times$, so we can choose distinct $r,s\in\{2,3,4\}$ with $r-s=a/b$. Then the weighted sum is $-a+b(r-s)=0$.

Otherwise Lemma~\ref{lem:four-weights} leaves four patterns. Table~\ref{tab:balance} gives a split and a coloring for every admissible vertex-weight pattern. Within a row, subclasses may be reordered to match the displayed weights. All dot products are divisible by five.
\begin{table}[ht]
\centering\small
\caption{Explicit independent-subclass assignments after one split. All weights and colors are interpreted modulo five; dot products are integer representatives.}\label{tab:balance}
\begin{tabular}{@{}cccc@{}}
\toprule
Original class weights & Split subclass weights & Colors & Dot product\\
\midrule
$(1,1,1,2)$ & $(1,1,1,1,1)$ & $(0,1,2,3,4)$ & $10$\\
$(1,3,3,3)$ & $(3,3,3,3,3)$ & $(0,1,2,3,4)$ & $30$\\
                  & $(1,1,2,3,3)$ & $(0,1,3,2,4)$ & $25$\\
$(2,2,2,4)$ & $(1,1,2,2,4)$ & $(0,1,2,4,3)$ & $25$\\
                  & $(2,2,3,4,4)$ & $(0,1,3,2,4)$ & $35$\\
$(3,4,4,4)$ & $(1,3,3,4,4)$ & $(0,1,4,2,3)$ & $35$\\
\bottomrule
\end{tabular}
\end{table}

To verify that a row always applies, observe first that in $(1,1,1,2)$ a weight-two class splits into two weight-one subclasses. For $(1,3,3,3)$, a nonsingleton weight-one class has pattern $(3,3)$ and yields the first row for this multiset. Otherwise $|V(J)|\ge5$ forces a nonsingleton weight-three class, whose pattern is $(1,1,1)$; this splits into weights one and two. For $(2,2,2,4)$, a weight-two class has pattern $(1,1)$ or $(1,3,3)$, allowing splits of weights $(1,1)$ or $(3,4)$. Finally, for $(3,4,4,4)$, any weight-four class splits into weights one and three. The displayed rows exhaust all cases.
\end{proof}

\begin{lemma}[Balancing modulo four]\label{lem:balance-four}
Let $J$ be four-colorable and assign weights $w(v)\in\Z/4\Z$. If $W=\sum_vw(v)\ne0\pmod4$, then $J$ admits a proper coloring $h:V(J)\to\Z/4\Z$ with $\sum_vw(v)h(v)=0\pmod4$.
\end{lemma}
\begin{proof}
Start with a proper four-coloring, allowing empty classes. If $W$ is odd, translating all four colors by a common residue changes the weighted sum by $Wt$, which attains every residue modulo four.

If $W=2$, let $W_0,\ldots,W_3$ be the class weights. Two of these weights have the same parity. Give those two classes the odd colors $1,3$ and the other classes the even colors $0,2$. The weighted sum $L$ is then even. A translation by $t$ replaces it with $L+2t$, so a suitable $t$ makes it zero. Both permutations and translations of auxiliary colors preserve properness.
\end{proof}

\begin{theorem}\label{thm:cubic}
Every cubic 3-connected plane graph satisfies $\lambda(G)\le14$. If $G$ is bipartite, then $\lambda(G)\le9$. Moreover, $\lambda(G)\le11$ in the one-class case of Lemma~\ref{lem:cokernel}, and also in its two-class case whenever
\begin{equation}\label{eq:modfour-condition}
|D|-2|A|\not\equiv0\pmod4,
\end{equation}
where the classes $A,D$ occur with multiplicities one and two in each dual triangle.
\end{theorem}
\begin{proof}
For $G=K_4$, $\lambda(G)=3$. Otherwise the connected cubic graph $G$ is not complete, so Brooks' theorem gives a proper three-coloring of $G$.

\emph{One class.}
If $s=1$, the incidence map is surjective modulo four because three is a unit in $\Z/4\Z$. Choose a proper four-coloring $h$ of $J$ by the Four Color Theorem~\cite{fourcolor}, and solve $Bb=h\pmod4$ using Lemma~\ref{lem:cokernel}. Lemma~\ref{lem:modular-combination} yields $\lambda(G)\le 3\cdot4-1=11$.

\emph{Two classes.}
If \eqref{eq:modfour-condition} holds, give dual vertices in $A$ weight $-2$ and those in $D$ weight $1$ modulo four. The total weight is nonzero, so Lemma~\ref{lem:balance-four} gives a proper four-coloring with
\[
\sum_{v\in D}h(v)=2\sum_{v\in A}h(v)\pmod4.
\]
Lemma~\ref{lem:cokernel} supplies $Bb=h\pmod4$, giving again $\lambda(G)\le11$.

In the remaining two-class case the total dual weight vanishes modulo four, so translation alone no longer guarantees balancing. Instead give vertices of $A$ weight $3=-2\pmod5$ and vertices of $D$ weight $1$. The dual has at least five vertices. By Lemma~\ref{lem:balance-five}, it admits a proper $h:V(J)\to\FF_5$ with $\sum_Dh=2\sum_Ah$. Solve $Bb=h\pmod5$ by Lemma~\ref{lem:cokernel}, then use Lemma~\ref{lem:modular-combination} to obtain $\lambda(G)\le3\cdot5-1=14$.

\emph{Three classes.}
If $s=3$, the equivalence classes $A,D,E$ properly three-color $J$. Each has at least two vertices. Indeed, if some class consisted of a single vertex $x$, every facial triangle would contain $x$. Then $2|V(J)|-4\le d_J(x)\le |V(J)|-1$, impossible for $|V(J)|\ge4$.

If two classes had size two, every face would contain an edge between them. There are at most four such edges, each incident with two faces, so $2|V(J)|-4\le8$ and $|V(J)|\le6$. Since each class has size at least two, all three have size two; $J$ is the octahedron and $G$ is the cube. For the cube, list its vertices in binary lexicographic order and assign the labels
\[
(0,1,1,3,2,0,3,2).
\]
They are proper, and the sums on the three pairs of opposite faces are $(5,7)$, $(3,9)$, and $(6,6)$. The equal sums of six belong to opposite, not edge-adjacent, faces; the three pairs use disjoint sets of values. Hence all edge-adjacent faces have distinct sums, proving $\lambda(G)\le3$ in this case.

Otherwise rename the classes so $|D|,|E|\ge3$. Put $h=0$ on $A$, assign values from $\{1,2\}$ on $D$, and values from $\{3,4\}$ on $E$. For each $m\ge3$, the intervals $[m,2m]$ and $[3m,4m]$ contain multiples of five: this is checked directly for $m=3,4$, and for $m\ge5$ follows from their lengths. Every integer in each interval is realizable as a sum of the specified values, so choose both class sums to vanish modulo five. The palettes are disjoint; hence $h$ is proper, and Lemma~\ref{lem:cokernel} gives $Bb=h\pmod5$.

A properly three-colored plane triangulation has even degrees, because the other two colors alternate along each vertex link. Therefore all faces of $G$ have even length and $G$ is bipartite. Apply Lemma~\ref{lem:modular-combination} with $r=2,q=5$ to conclude $\lambda(G)\le9$.

Conversely, if $G$ is cubic and bipartite, its dual triangulation has even degrees and is properly three-colorable by the classical criterion for plane triangulations (see~\cite{diestel}). Such a coloring is constant on each opposite-vertex equivalence class, so the dual must be in the $s=3$ case. This proves the claimed bound for every bipartite cubic graph.
\end{proof}

\section{An eight-element list theorem for fullerenes}\label{sec:fullerenes}
\subsection{Local face structure}
\begin{definition}
A \emph{fullerene} is a 3-connected cubic plane graph all of whose faces have length five or six.
\end{definition}
For a fullerene with $n$ vertices and $f_5,f_6$ faces of each length, Euler's formula and cubicity give
\[
f_5+f_6=\frac n2+2,\qquad5f_5+6f_6=3n,
\qquad f_5=12,\qquad f_6=\frac n2-10.
\]
In particular, $n\ge20$, and a pentagonal face exists. For fullerene generation and computational background, see Brinkmann, Goedgebeur, and McKay~\cite{brinkmann}. We use the theorem that every fullerene is \emph{cyclically five-edge-connected}: no set of fewer than five edges separates two components both containing cycles. For a corrected and simplified proof, see Qi and Zhang~\cite{qi}, who repair an oversight in the earlier proof of Do\v{s}li\'c~\cite{doslic}.

\begin{lemma}\label{lem:pentagon-neighbors}
The dual of a fullerene has no separating triangle. If $C=v_0v_1v_2v_3v_4v_0$ is a pentagonal face and $H_i$ is the other face at $v_iv_{i+1}$ (indices modulo five), then $H_0,\ldots,H_4$ are distinct. Consecutive $H_i$ share the outgoing edge at their common vertex of $C$, and nonconsecutive $H_i$ share no edge.
\end{lemma}
\begin{proof}
A separating triangle of the dual would give a three-edge bond in the primal. If one side of this bond were a tree on $a$ vertices, cubicity would force $3a=2(a-1)+3$, hence $a=1$; the corresponding dual triangle would then be facial. Thus both sides of a \emph{separating} dual triangle would contain cycles, contrary to cyclic five-edge-connectivity.

The faces $H_i$ are distinct by the polyhedral face-intersection property, and consecutive adjacencies follow from cubicity. In the dual, $H_0,\ldots,H_4$ form the link cycle of the vertex corresponding to $C$. Any edge between nonconsecutive link vertices, together with that central vertex, would form a separating triangle, since both intervening link arcs have internal vertices. Such an edge is therefore impossible.
\end{proof}

\subsection{Greedy labeling outside a pentagon}
Index the two constraints at every original edge $e$ by $(e,\mathrm{edge})$ and $(e,\mathrm{face})$. Their linear factors are $x_p-x_q$ for $e=pq$ and $S_F-S_H$ for its incident faces. For a linear factor $L$, its \emph{support} consists of the variables whose coefficients are nonzero after cancellation. We call a factor occurrence \emph{determined} when every variable in its support has received a label. Repeated, identical linear expressions indexed by different constraints remain separate \emph{occurrences}; counting all of them can only overestimate forbidden choices.

For $e=pq$, the two incident faces intersect exactly in the edge $pq$. Hence the endpoints cancel from $S_F-S_H$, while each other boundary vertex occurs with coefficient $+1$ or $-1$. Once all but one supported variable $x_u$ are labeled, such a factor becomes a nonconstant affine function of $x_u$ and excludes at most one real value.

\begin{lemma}\label{lem:outside-pentagon}
Fix a pentagonal face $C$ of a fullerene. From arbitrary eight-element real lists one can label all vertices outside $C$ so that every determined constraint is satisfied.
\end{lemma}
\begin{proof}
Leave the five vertices of $C$ unlabeled, and process the remaining vertices in nonincreasing distance from $C$, breaking ties arbitrarily. Maintain the invariant that every determined factor is nonzero.

When processing $u$, choose a neighbor $v$ at strictly smaller distance; $v$ is still unlabeled. If $v\notin C$, choose a neighbor $w$ of $v$ at still smaller distance. If $v\in C$, choose either neighbor of $v$ along $C$ as $w$. In either case $u,v,w$ are distinct and $v,w$ remain unlabeled.

Let $A,B$ be the two faces at $uv$, with lengths $a,b$, and let $R$ be the third face at $u$, with length $c$. Three edge-factor occurrences contain $x_u$. For each of $A,B,R$, the facial-factor occurrences containing $x_u$ are indexed by the boundary edges not incident with $u$. These occurrences are distinct: different faces containing $u$ share only the edge incident with $u$, where $x_u$ cancels. The total number of occurrences supported on $u$ is therefore
\[
M_u=3+(a-2)+(b-2)+(c-2)=a+b+c-3.
\]
We now identify occurrences that remain \emph{undetermined} after $u$ is labeled.

\emph{First family: $a+b-5$ occurrences containing $x_v$.}
The edge factor at $uv$ remains undetermined. So do the $a-3$ facial factors indexed by edges of $A$ incident with neither $u$ nor $v$, and the analogous $b-3$ factors for $B$. Their opposite faces contain neither $u$ nor $v$ by the polyhedral intersection property; thus these facial differences retain $x_v$ with a nonzero coefficient. These $1+(a-3)+(b-3)=a+b-5$ distinct occurrences all contain $x_v$ in their supports.

\emph{Second family: one further occurrence containing $x_w$ but not $x_v$.}
Let $t$ be the third neighbor of $v$, other than $u,w$. The two faces $P,Q$ at $vt$ contain $u,w$, respectively, and intersect exactly in $vt$. Thus $S_P-S_Q$ has coefficients $+1,-1$ on $x_u,x_w$, up to an overall sign, but coefficient zero on $x_v$. It remains undetermined because $w$ is unlabeled. It is distinct from the first family because every occurrence there has nonzero $x_v$-coefficient.

Only factor occurrences supported on $u$ can become newly determined when we label $u$. At least $(a+b-5)+1=a+b-4$ of the $M_u$ occurrences are still undetermined. Hence at most
\[
M_u-(a+b-5)-1=c+1\le7
\]
new constraints can become determined. Each excludes at most one value of $x_u$. The eight-element list at $u$ therefore contains a permissible choice, preserving the invariant. Iterating completes the labeling outside $C$.
\end{proof}

\subsection{A uniform polynomial extension across the pentagon}
\begin{figure}[ht]
\centering
\begin{tikzpicture}[
  scale=1.12,
  font=\small,
  core/.style={circle,fill=black,inner sep=2pt},
  outside/.style={
    circle,draw=black,fill=white,
    line width=0.8pt,inner sep=2pt
  }
]

  \foreach \i in {0,...,4} {
    \pgfmathsetmacro{\a}{90-72*\i}
    \coordinate (v\i) at (\a:1.52);
    \coordinate (t\i) at (\a:2.38);
  }

  \draw[thick]
    (v0)--(v1)--(v2)--(v3)--(v4)--cycle;

  \foreach \i in {0,...,4} {
    \draw[thick] (v\i)--(t\i);
  }

  \foreach \i in {0,...,4} {
    \node[core] at (v\i) {};
    \node[outside] at (t\i) {};
  }

  \node[left=5pt]        at (v0) {$v_0$};
  \node[below right=4pt] at (v1) {$v_1$};
  \node[right=5pt]       at (v2) {$v_2$};
  \node[left=5pt]        at (v3) {$v_3$};
  \node[below left=4pt]  at (v4) {$v_4$};

  \node[above=5pt]       at (t0) {$t_0$};
  \node[above right=4pt] at (t1) {$t_1$};
  \node[below right=4pt] at (t2) {$t_2$};
  \node[below left=4pt]  at (t3) {$t_3$};
  \node[above left=4pt]  at (t4) {$t_4$};

  \node[font=\large] at (0,0) {$C$};

  \foreach \i in {0,...,4} {
    \pgfmathsetmacro{\a}{90-72*\i-36}
    \node at (\a:2.05) {$H_{\i}$};
  }

  \node[
    align=left,
    text width=4.1cm,
    anchor=west,
    font=\footnotesize
  ] at (3.3,0.2) {
    The two faces at $v_iv_{i+1}$ are
    $C$ and $H_i$.\\[5pt]
    Each $H_i$ meets $C$ exactly in
    $\{v_i,v_{i+1}\}$.\\[5pt]
    All labels at open vertices are fixed
    before the extension.
  };
\end{tikzpicture}
\caption{Schematic neighborhood of the unlabeled pentagonal face. Filled vertices carry the five remaining variables; open vertices have already been labeled. Each surrounding $H_i$ has its edge on $C$, two outgoing edges, and an outer boundary path of length $2+\varepsilon_i$. Outer face boundaries are omitted; comparisons across their edges yield the final row of Table~\ref{tab:pentagon}.}\label{fig:pentagon}
\end{figure}
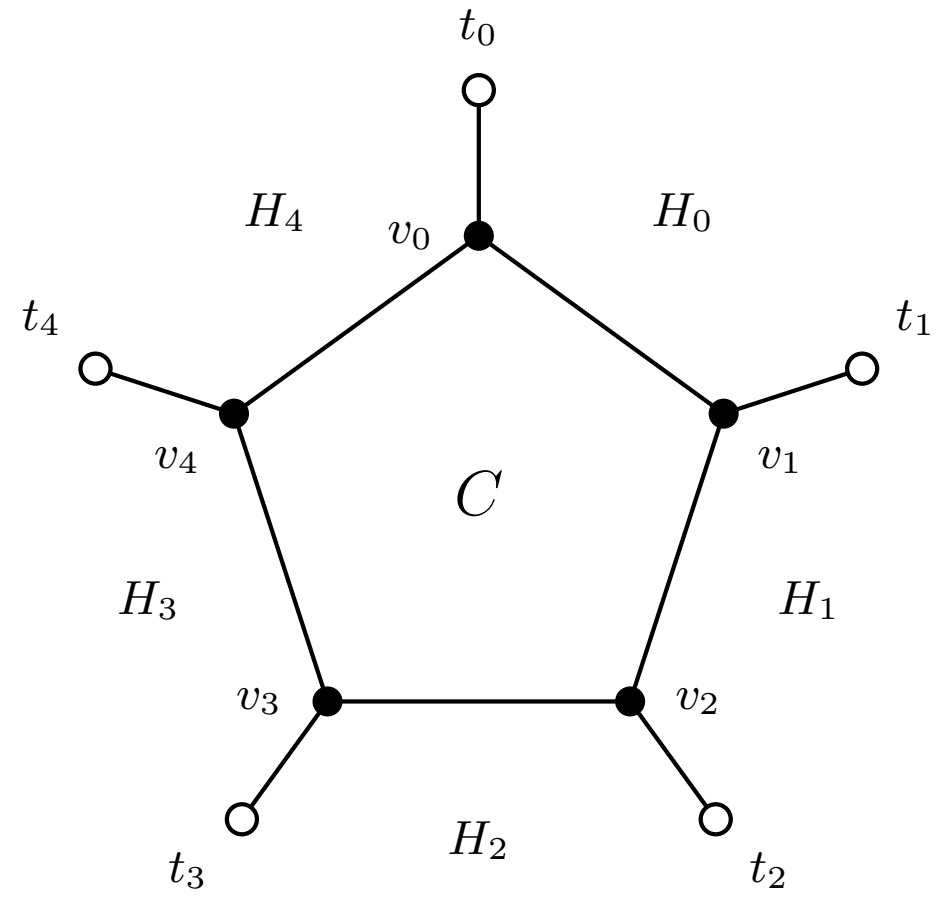

The remaining constraints involve five variables and a number of fixed constants determined by the preceding greedy stage. Figure~\ref{fig:pentagon} identifies the relevant faces and illustrates the geometry.

\begin{lemma}\label{lem:pentagon-extension}
Any labeling of $V(G)\setminus V(C)$ satisfying every determined constraint extends to $C$ from arbitrary eight-element real lists.
\end{lemma}
\begin{proof}
Write $x_i=c(v_i)$, with indices modulo five, and let $|H_i|=5+\varepsilon_i$, $\varepsilon_i\in\{0,1\}$. Let $t_i$ be the outside neighbor of $v_i$, put $a_i=c(t_i)$, and let $b_i$ be the sum of labels on $V(H_i)\setminus V(C)$. Finally, set $S=x_0+\cdots+x_4$.

By Lemma~\ref{lem:pentagon-neighbors}, all remaining constraints are listed in Table~\ref{tab:pentagon}. The table records \emph{factor occurrences}, so two different edges are still counted separately even when they yield the same expression after some labels have been fixed.
\begin{table}[ht]
\centering\small
\caption{All constraint occurrences that involve a vertex of $C$. Each affine form is specified up to sign; indices are modulo five.}\label{tab:pentagon}
\begin{tabular}{@{}>{\raggedright\arraybackslash}p{3.05cm}c p{7.0cm}@{}}
\toprule
Constraint family & Number & Affine form\\
\midrule
Outgoing edges & $5$ & $x_i-a_i$\\
Edges of $C$ & $5$ & $x_i-x_{i+1}$\\
Comparisons $H_{i-1},H_i$ & $5$ & $x_{i-1}-x_{i+1}+b_{i-1}-b_i$\\
Comparisons $C,H_i$ & $5$ & $S-x_i-x_{i+1}-b_i$\\
Other comparisons at $H_i$ & $\sum_i(2+\varepsilon_i)$ & $x_i+x_{i+1}+d_{ij}$, with fixed $d_{ij}$\\
\bottomrule
\end{tabular}
\end{table}

Here is an edge-by-edge check of completeness. The only edges incident with $C$ are the five cycle edges and five outgoing edges $v_it_i$, giving the first two rows. At each cubic vertex $v_i$, the incident faces are exactly $C,H_{i-1},H_i$. The cycle edge $v_iv_{i+1}$ gives the comparison $C,H_i$, and the outgoing edge $v_it_i$ gives $H_{i-1},H_i$, producing the next two rows.

The boundary of $H_i$ consists of its edge on $C$, the two outgoing edges at $v_i,v_{i+1}$, and $|H_i|-3=2+\varepsilon_i$ other edges. If $e$ is one of these other edges and $K$ is the opposite face, then $K\ne C$. Moreover, $K$ has no vertex of $C$: if $K$ contained some $v_j$, cubicity would give $K=H_{j-1}$ or $H_j$; Lemma~\ref{lem:pentagon-neighbors} rules out nonconsecutive surrounding faces sharing an edge, while consecutive surrounding faces share only their already counted outgoing edge. As $V(H_i)\cap V(C)=\{v_i,v_{i+1}\}$, this comparison is exactly
\[
S_{H_i}-S_K=x_i+x_{i+1}+b_i-s_c(K),
\]
with $d_{ij}=b_i-s_c(K)$ fixed. No other edge contributes a new factor involving $C$, and an other edge cannot be counted for two distinct $H_i$. All remaining constraints were already satisfied outside the pentagon. In particular, each displayed affine factor has a nonzero coefficient on at least one $x_i$.

Let $P(x)$ be the product of these factors. Its degree is $30+\sum_i\varepsilon_i$. The two groups of difference factors supply, in the highest-degree homogeneous part, all ten pairwise differences $x_i-x_j$ for $0\le i<j\le4$. Thus the top-degree part of $P$, up to sign, is
\[
\biggl(\prod_{i=0}^4x_i\biggr)
\biggl(\prod_{0\le i<j\le4}(x_i-x_j)\biggr)
\biggl(\prod_{i=0}^4(S-x_i-x_{i+1})\biggr)
\biggl(\prod_{i=0}^4(x_i+x_{i+1})^{2+\varepsilon_i}\biggr).
\]
The differences between nonconsecutive indices arise from facial comparisons; they do not impose vertex properness between nonadjacent vertices of $C$.

For each $i$ with $\varepsilon_i=0$, multiply $P$ by the auxiliary factor $x_i+x_{i+1}$. The resulting polynomial has degree $35$ and top homogeneous part
\begin{equation}\label{eq:Qtilde}
\widetilde Q(x)=
\biggl(\prod_{i=0}^4x_i\biggr)
\biggl(\prod_{0\le i<j\le4}(x_i-x_j)\biggr)
\biggl(\prod_{i=0}^4(S-x_i-x_{i+1})\biggr)
\biggl(\prod_{i=0}^4(x_i+x_{i+1})^3\biggr),
\end{equation}
up to an overall sign. Appendix~\ref{app:coefficient} proves the symbolic identity
\begin{equation}\label{eq:coefficient}
[x_0^7x_1^7x_2^7x_3^7x_4^7]\,\widetilde Q\equiv2\pmod7.
\end{equation}
This integer coefficient is nonzero over $\R$. The Combinatorial Nullstellensatz~\cite{alon}, with exponent seven in each variable, therefore guarantees a nonzero evaluation on the product of the five prescribed eight-element real lists. The extra factors impose only additional nonvanishing conditions, so nonvanishing of the augmented polynomial forces every original factor of $P$ to be nonzero. This completes the extension.
\end{proof}

\begin{theorem}\label{thm:fullerene-list}
Every fullerene $G$ satisfies $\laml(G)\le7$, and consequently $\lambda(G)\le7$.
\end{theorem}
\begin{proof}
First apply Lemma~\ref{lem:outside-pentagon} outside a pentagonal face, then Lemma~\ref{lem:pentagon-extension} to its boundary vertices.
\end{proof}

\subsection{The dodecahedron and the common-alphabet question}
\begin{proposition}\label{prop:dodecahedron}
For the dodecahedral graph $D$, $\lambda(D)=3$. Consequently, $3\le K_{\Ful}\le7$.
\end{proposition}
\begin{proof}
For a proper labeling of a five-cycle from $\{0,1,2\}$, the label multiplicities are $2,2,1$. If the singleton label is $a$, the face sum is $6-a$. Thus distinct sums on edge-adjacent pentagonal faces would give a proper three-coloring of the dual icosahedron. But the five neighbors of any icosahedral vertex induce an odd cycle and would have to use only the other two colors. This is impossible, so $\lambda(D)\ge3$.

For the upper bound, use $D=GP(10,2)$ with edges $u_iu_{i+1}$, $u_iv_i$, and $v_iv_{i+2}$ (indices modulo ten). Assign the labels in Table~\ref{tab:dodecahedron}.
\begin{center}
\begin{minipage}{\textwidth}
\centering\small
\captionsetup{hypcap=false}
\captionof{table}{A four-label solution for the dodecahedron. Here $F_i=u_iu_{i+1}u_{i+2}v_{i+2}v_i$.}\label{tab:dodecahedron}
\begin{tabular}{@{}c*{10}{c}@{}}
\toprule
$i$&0&1&2&3&4&5&6&7&8&9\\\midrule
$c(u_i)$&2&3&2&0&2&3&2&1&0&3\\
$c(v_i)$&0&0&1&2&0&1&3&0&2&2\\
$s_c(F_i)$&8&7&5&8&10&7&8&6&7&10\\\bottomrule
\end{tabular}
\end{minipage}
\end{center}
The two remaining faces are the even and odd inner pentagons, of sums six and five, respectively. All edges of the three specified types are properly labeled. Each $F_i$ is adjacent to $F_{i\pm1}$, $F_{i\pm2}$, and the inner pentagon of the same parity. The table gives distinct sums on every such adjacency. Therefore $\lambda(D)\le3$.
\end{proof}

\begin{conjecture}[Common-alphabet fullerene problem]\label{conj:fullerene}
Every fullerene satisfies $\lambda(G)\le3$, equivalently $K_{\Ful}=3$.
\end{conjecture}

\begin{remark}[The distinct list problem]\label{rem:fullerene-list}
Can the bound $\sup_{G\in\Ful}\laml(G)\le7$ be reduced, for example to six? This question concerns arbitrary real lists, whereas Conjecture~\ref{conj:fullerene} permits only the common alphabet $\{0,1,2,3\}$. Proving the common-alphabet conjecture alone would not improve the list bound, although a stronger arbitrary-list theorem would imply its corresponding common-alphabet statement. The pentagon extension and coefficient method might admit further compression for a common short integer alphabet, but no such improvement is established here.
\end{remark}

\section{A universal constant for polyhedral plane graphs}\label{sec:universal}
For each edge $e$, order its incident faces as $F_e^+,F_e^-$ and define the integer facial-difference vector
\begin{equation}\label{eq:rho}
\rho_e(v)=\ind_{\{v\in V(F_e^+)\}}-\ind_{\{v\in V(F_e^-)\}}.
\end{equation}
Each $\rho_e$ is a nonzero vector in $\{-1,0,1\}^{V(G)}$, and the corresponding facial condition is $\langle\rho_e,c\rangle\ne0$. We first bound the number of distinct projections of these vectors, identifying a vector with its negative. The proof combines two plane constructions arising from the same facial types: a skeleton on selected primal vertices and a quotient of an incidence--dual auxiliary graph.

\subsection{Counting projected facial-difference types}
Fix a nonempty $U\subseteq V(G)$ and set $r=|U|$. The \emph{type} of a face $F$ is
\[
A_F=V(F)\cap U.
\]
Let $N(U)$ count the distinct nonzero vectors $\ind_{A_F}-\ind_{A_H}$ arising from edge-adjacent faces, with opposite vectors identified. For $r\ge3$, let $t_2,t_3,h$ be the numbers of \emph{distinct types} of sizes two, three, and at least four. Put
\[
s=\sum_{F:\,|A_F|\ge3}|A_F|.
\]
Thus $s$ counts selected-vertex incidences in original faces of type size at least three, not simply the number of distinct types. Let $z$ count the original faces $F$ with $|A_F|\ge3$ that have an edge-adjacent face of empty type. Let $P$ count the unordered selected pairs that occur either as a two-element type or as a subset of a three-element type.

Two distinct faces containing the same pair of vertices must share their original connecting edge, by the polyhedral face-intersection property. Hence any two-element type occurs in at most two faces, and if twice, those faces are edge-adjacent. A type of size at least three occurs in only one face. These facts will justify the counting multiplicities.

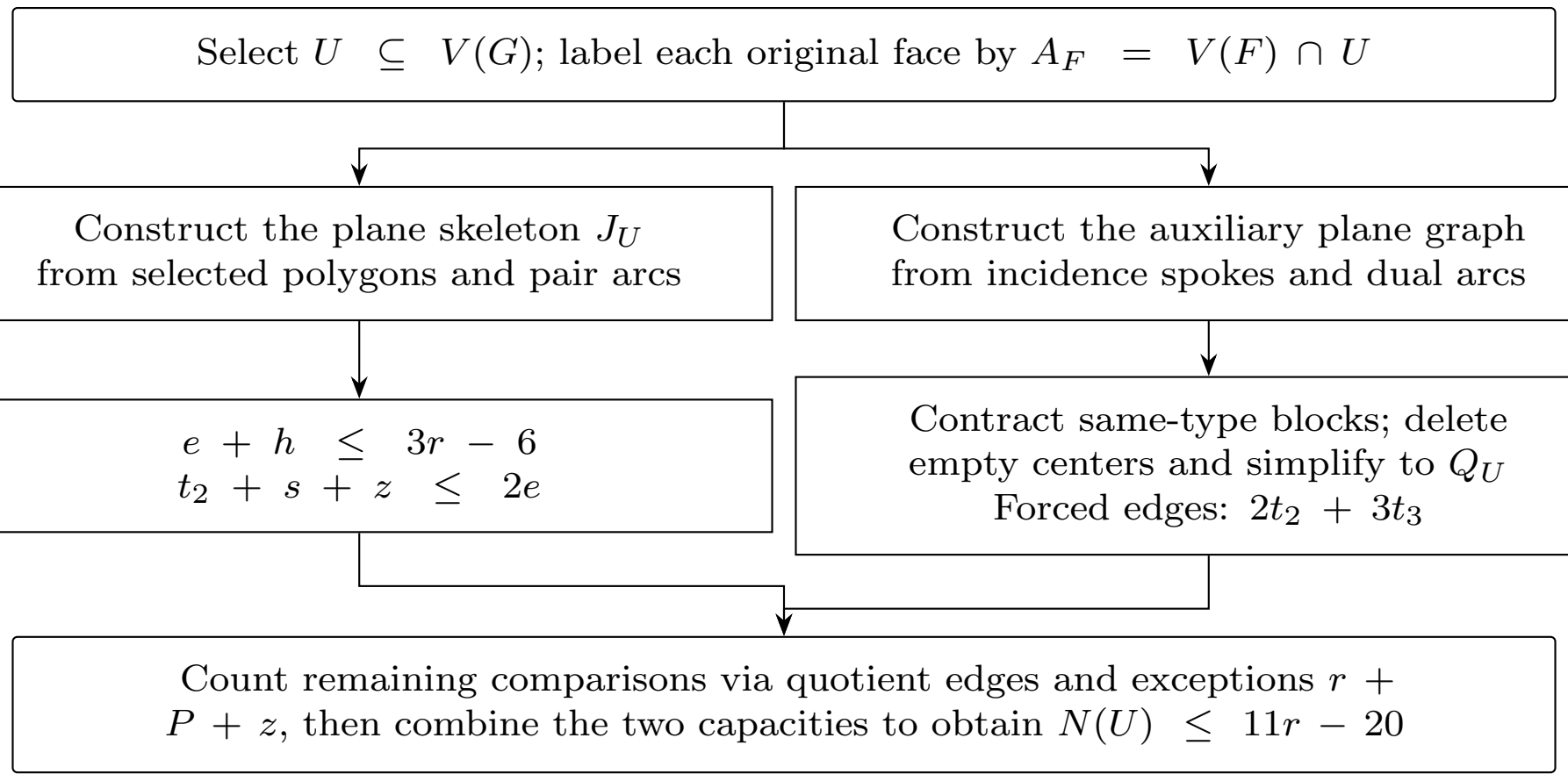
\begin{figure}[ht]
\centering
\begin{tikzpicture}[font=\footnotesize,>=Stealth,box/.style={draw,rounded corners=1pt,align=center,text width=5.8cm,inner sep=6pt},stepbox/.style={draw,align=center,text width=5.8cm,inner sep=6pt}]
\node[box,text width=11.2cm] (start) at (0,0) {Select $U\subseteq V(G)$; label each original face by $A_F=V(F)\cap U$};
\node[stepbox] (left) at (-3.2,-1.5) {Construct the plane skeleton $J_U$\\from selected polygons and pair arcs};
\node[stepbox] (right) at (3.2,-1.5) {Construct the auxiliary plane graph\\from incidence spokes and dual arcs};
\draw[->] (start.south) -- ++(0,-0.35) -| (left.north);
\draw[->] (start.south) -- ++(0,-0.35) -| (right.north);
\node[stepbox] (lout) at (-3.2,-3.1) {$e+h\le3r-6$\\ $t_2+s+z\le2e$};
\node[stepbox] (rout) at (3.2,-3.1) {Contract same-type blocks; delete\\empty centers and simplify to $Q_U$\\Forced edges: $2t_2+3t_3$};
\draw[->] (left)--(lout);
\draw[->] (right)--(rout);
\node[box,text width=11.2cm] (last) at (0,-4.9) {Count remaining comparisons via quotient edges and exceptions $r+P+z$, then combine the two capacities to obtain $N(U)\le11r-20$};
\draw[->] (lout.south) -- ++(0,-0.4) -| (last.north);
\draw[->] (rout.south) -- ++(0,-0.4) -| (last.north);
\end{tikzpicture}
\caption{Proof architecture for the projected-type estimate. The left branch controls available incidence capacity; the right branch counts nonexceptional comparison vectors using a simple plane quotient. Arrows indicate logical dependencies, not additional graph contractions.}\label{fig:architecture}
\end{figure}

\begin{lemma}[Skeleton and unused incidence capacity]\label{lem:skeleton}
For $r\ge3$ there is a simple plane graph $J_U$ on $U$, with $e=|E(J_U)|$, satisfying
\begin{align}
e+h&\le3r-6,\label{eq:skeleton1}\\
t_2+s+z&\le2e,\qquad t_2+3t_3+4h+z\le2e,\label{eq:skeleton2}\\
t_2&\le P\le e,\qquad z\le t_3+h.\label{eq:skeleton3}
\end{align}
\end{lemma}
\begin{proof}
Within each original face containing at least three selected vertices, draw a polygon joining those vertices in cyclic boundary order. Its sides are arcs through the face interior, except that a side joining originally consecutive vertices may be moved onto their boundary edge. For every distinct two-element type, choose one representative face and draw an arc connecting its two selected vertices. Arc interiors are disjoint within each original face. If arcs from distinct faces have the same endpoints $u,v$, those faces share the original edge $uv$; move both arcs onto that edge and identify them. Since no third face contains $u,v$, the resulting graph $J_U$ is simple and plane.

A selected polygon with at least four vertices bounds a region containing no edge of $J_U$ in its interior. Choose one diagonal between nonconsecutive selected vertices in each such polygon. It cannot duplicate an existing edge of $J_U$: it is not a side of that polygon, and any other original face containing both endpoints would have to share their original boundary edge, contrary to nonconsecutiveness. Diagonals from distinct polygons are also distinct, by the same face-intersection property, and can be drawn without crossings. Adding all $h$ diagonals yields a simple plane graph with $e+h$ edges on $r\ge3$ vertices. The planar edge bound gives \eqref{eq:skeleton1}, without requiring connectedness.

For the capacity count, a \emph{contributor} is either a pair consisting of an original face and one side of its selected polygon, or the single chosen arc for a two-element type. Each contributor is assigned to the edge of $J_U$ with its endpoints. There are exactly $t_2+s$ contributors. A fixed original face contributes at most once to a given endpoint pair, and at most two faces can contain that pair. Consequently the number $c_f$ of contributors assigned to $f\in E(J_U)$ is either one or two, and
\[
\sum_{f\in E(J_U)}c_f=t_2+s,\qquad
\sum_{f\in E(J_U)}(2-c_f)=2e-(t_2+s).
\]

We inject the $z$ marked faces into the edges with one unused contributor slot. For each marked face $F$, choose a boundary edge $xy$ shared with a face of empty type. Both $x,y$ lie outside $U$. Let $u,v$ be the consecutive selected vertices that bracket the boundary stretch of $F$ containing $xy$. That stretch has at least three original edges. The complementary boundary stretch contains another selected vertex, because $|A_F|\ge3$. Therefore $uv$ is not an original boundary edge of $F$. No second face can contain both $u,v$: by the polyhedral intersection property it would have to share the original edge $uv$ with $F$, which does not exist on the boundary of $F$. The selected polygon side $uv$ thus has exactly one contributor, so $c_{uv}=1$.

The assignment $F\mapsto uv$ is injective. Two distinct marked faces cannot select the same pair, because that pair occurs in no other original face. Hence
\[
z\le\sum_{f\in E(J_U)}(2-c_f),
\]
which gives $t_2+s+z\le2e$. Since $s\ge3t_3+4h$ (each type of size at least three occurs only once), the second inequality of \eqref{eq:skeleton2} follows.

Every pair counted by $P$ appears as a side of a chosen two-element arc or a three-vertex polygon, so $t_2\le P\le e$. Finally, every marked face has its own type of size at least three, giving $z\le t_3+h$. This proves \eqref{eq:skeleton3}.
\end{proof}

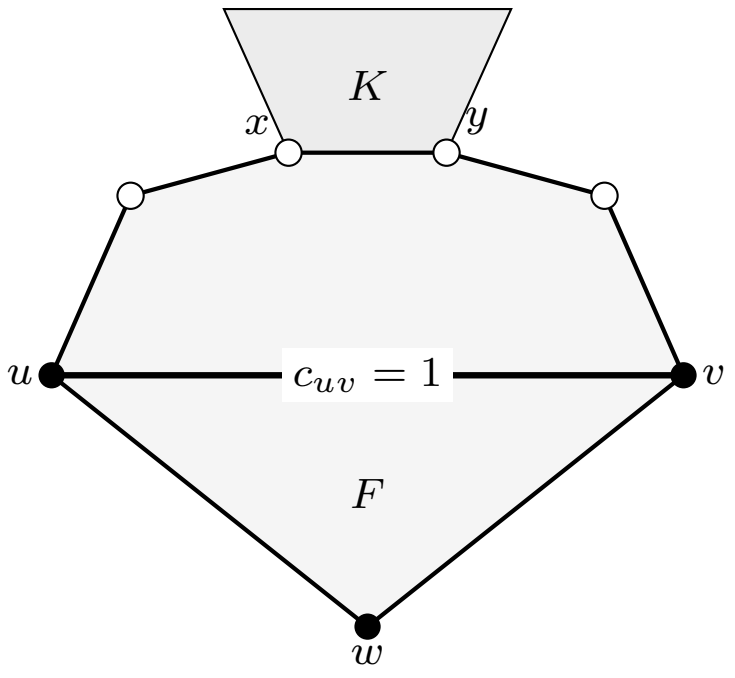
\begin{figure}[ht]
\centering
\begin{tikzpicture}[font=\footnotesize,scale=0.95]
  \coordinate (u) at (-2.2,0.1);
  \coordinate (p) at (-1.65,1.35);
  \coordinate (x) at (-0.55,1.65);
  \coordinate (y) at (0.55,1.65);
  \coordinate (q) at (1.65,1.35);
  \coordinate (v) at (2.2,0.1);
  \coordinate (w) at (0,-1.65);
  \fill[gray!8] (u)--(p)--(x)--(y)--(q)--(v)--(w)--cycle;
  \fill[gray!15] (x)--(-1.0,2.65)--(1.0,2.65)--(y)--cycle;
  \draw[thick] (u)--(p)--(x)--(y)--(q)--(v)--(w)--cycle;
  \draw[thin] (x)--(-1.0,2.65)--(1.0,2.65)--(y);
  \draw[very thick] (u)--(v);
  \draw[densely dashed] (u)--(w)--(v);
  \foreach \pt in {u,v,w}{\fill (\pt) circle (2.6pt);}
  \foreach \pt in {p,x,y,q}{\draw[fill=white] (\pt) circle (2.6pt);}
  \node[left] at (u) {$u$};\node[right] at (v) {$v$};
  \node[below] at (w) {$w$};
  \node[above left] at (x) {$x$};\node[above right] at (y) {$y$};
  \node at (0,-0.72) {$F$};
  \node at (0,2.12) {$K$};
  \node[fill=white,inner sep=2pt] at (0,0.10) {$c_{uv}=1$};
  \node[align=left,text width=4.2cm,anchor=west] at (3.05,0.5)
    {$xy$ borders an empty-type face $K$ ($A_K=\varnothing$).\\[5pt]
     In $F$, the selected set is $A_F=\{u,v,w\}$. The polygon side $uv$ joins consecutive \emph{selected} vertices, but is not an original boundary edge of $F$.\\[5pt]
     No other face can contribute to $uv$.};
\end{tikzpicture}
\caption{The unused-slot injection in Lemma~\ref{lem:skeleton}. Black vertices are selected and white vertices are not. The bold side $uv$ has exactly one contributor, from the marked face $F$, because it is not an original edge. Distinct marked faces must choose distinct such sides. The surrounding geometry is schematic.}\label{fig:unused-slot}
\end{figure}

\begin{figure}[ht]
\centering
\begin{tikzpicture}[font=\scriptsize,scale=0.97,
  fac/.style={circle,draw,fill=white,inner sep=2pt},
  vert/.style={circle,fill=black,inner sep=2pt}]
\begin{scope}[xshift=0cm]
 \node[fac,label=above:{$f$}] (f) at (0,1.25) {};
 \node[fac,label=below:{$g$}] (g) at (0,0) {};
 \node[vert,label=below:{$u$}] (u) at (1.15,-0.15) {};
 \node[vert,label=above:{$v$}] (v) at (1.25,1.95) {};
 \node[vert,label=left:{$w$}] (w) at (-1.1,1.6) {};
 \draw (f)--(v) (f)--(w) (f)--(u) (f)--(g) (g)--(u);
 \node[align=center] at (0,-1.05) {(a) Before contraction\\$A_f=\{u,v,w\}$, $A_g=\{u\}$};
\end{scope}
\draw[-{Stealth[length=2mm]}] (2,-0.1) -- (2.85,-0.1) node[midway,above] {$ug$};
\begin{scope}[xshift=4.0cm]
 \node[fac,label=above:{$f$}] (f) at (0,1.25) {};
 \node[vert,label=below:{$u=g$}] (u) at (1.15,-0.15) {};
 \node[vert,label=above:{$v$}] (v) at (1.25,1.95) {};
 \node[vert,label=left:{$w$}] (w) at (-1.1,1.6) {};
 \draw (f)--(v) (f)--(w);
 \draw (f) to[bend left=14] (u); \draw (f) to[bend right=14] (u);
 \node[align=center] at (0,-1.05) {(b) Multigraph\\two parallel $fu$ edges};
\end{scope}
\draw[-{Stealth[length=2mm]}] (6.05,-0.1) -- (6.85,-0.1) node[midway,above] {merge};
\begin{scope}[xshift=8.0cm]
 \node[fac,label=above:{$f$}] (f) at (0,1.25) {};
 \node[vert,label=below:{$u$}] (u) at (1.15,-0.15) {};
 \node[vert,label=above:{$v$}] (v) at (1.25,1.95) {};
 \node[vert,label=left:{$w$}] (w) at (-1.1,1.6) {};
 \draw[thick] (f)--(v) (f)--(w) (f)--(u);
 \node[align=center] at (0,-1.05) {(c) Simple quotient\\$\ind_{A_f}-e_u=e_v+e_w$};
\end{scope}
\end{tikzpicture}
\caption{Local simplification of a comparison $A_f=\{u,v,w\}$ versus $A_g=\{u\}$. Contracting the singleton incidence $ug$ makes the dual comparison parallel to the forced edge $fu$. The resulting pair indicator is counted as an exceptional type, so merging these edges loses no uncounted vector.}\label{fig:quotient-local}
\end{figure}
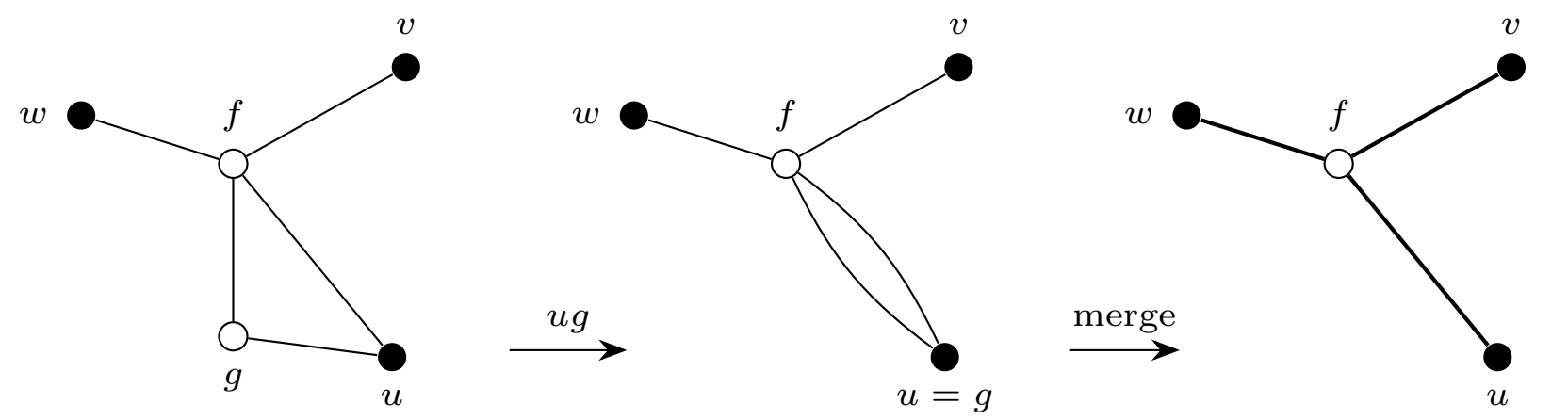

\begin{lemma}[A plane quotient with forced incidence edges]\label{lem:quotient}
For $r\ge3$ there is a simple plane graph $Q_U$ with
\[
R=r+t_2+t_3+h
\]
vertices. Its vertices have distinct nonempty type labels: each $u\in U$ has label $\{u\}$, and each distinct type $A$ with $|A|\ge2$ has one vertex $q_A$. The graph contains a set $E_0$ of $2t_2+3t_3$ forced incidence edges, and
\begin{equation}\label{eq:quotientbound}
|E(Q_U)\setminus E_0|\le3R-6-2t_2-3t_3.
\end{equation}
Every comparison between two nonempty types either has equal labels or is represented by an edge of $Q_U$ joining its two type labels.
\end{lemma}
\begin{proof}
Construct an auxiliary plane graph as follows. Place a center $c_F$ inside every face $F$. Join it by incidence spokes to each selected vertex in $V(F)\cap U$, and draw a dual edge $c_Fc_H$ through each original edge separating adjacent faces $F,H$. These arcs can be drawn simultaneously without crossings: within each face disk, arrange the spokes to selected boundary vertices and to an interior point of every boundary edge in the same cyclic order as their endpoints, then continue the latter spokes across the corresponding original edges. The original primal edges themselves are not edges of this auxiliary graph.

Now contract only connected blocks whose centers share the same type. The operations are specified in Table~\ref{tab:contractions}. For a singleton type $\{u\}$, contract the connected star
\[
K_u=\{u\}\cup\{c_F:A_F=\{u\}\}
\]
along its edges $uc_F$. Different such blocks are disjoint. If a two-element type occurs twice, the two corresponding faces are edge-adjacent, so contract the dual edge joining their centers to form one vertex $q_A$. These two-center blocks are disjoint from one another and from all singleton blocks. Retain the unique center of each type of size at least three and the center of a two-element type occurring only once. Delete all empty-type centers and their incident edges.

\begin{table}[ht]
\centering\footnotesize
\caption{Contraction and deletion rules. Each contraction block has a single type label.}\label{tab:contractions}
\begin{tabularx}{\textwidth}{@{}>{\raggedright\arraybackslash}p{3.1cm}>{\raggedright\arraybackslash}X>{\raggedright\arraybackslash}p{2.55cm}@{}}
\toprule
Original type & Operation & Surviving vertex\\
\midrule
Empty & Delete center and incident edges & None\\
$\{u\}$, any multiplicity & Contract all such centers with $u$ along the star $K_u$ & $u$ ($\{u\}$)\\
$\{u,v\}$, occurring twice & Contract the common dual edge of the two centers; later merge duplicate incidences & $q_{\{u,v\}}$\\
$|A|=2$ once, or $|A|\ge3$ & Retain the unique center & $q_A$\\
\bottomrule
\end{tabularx}
\end{table}

The empty-center deletions may be done first, since no empty center belongs to a prescribed block. Throughout the contractions, maintain the invariant that each surviving vertex represents exactly one nonempty type and that distinct types are never identified. Once every prescribed block is contracted, precisely one vertex remains for each nonempty type. Contraction of disjoint connected plane subgraphs and deletion preserve planarity; the result is a plane multigraph $Q_0$.

Delete loops and merge parallel edges to obtain a simple plane graph $Q_U$. A loop corresponds to two identical type labels and hence to the zero vector. All parallel comparison edges join the same two types, so their projected differences agree up to sign. Thus no distinct nonzero comparison between nonempty types is lost when parallel edges are merged. Figure~\ref{fig:quotient-local} illustrates a comparison edge becoming parallel to an incidence edge.

For each type $A$ of size two or three, every incidence edge $q_Au$ with $u\in A$ survives, possibly after merging duplicates when a two-element type occurs twice. These edges are not loops, and their pairs of endpoint labels are distinct. Mark all of them as $E_0$; there are exactly $2t_2+3t_3$. A comparison edge parallel to one of them does not add an edge to $Q_U$; its vector will be counted separately in Lemma~\ref{lem:encoding}. Since $R\ge r\ge3$, the simple planar edge bound $|E(Q_U)|\le3R-6$ proves \eqref{eq:quotientbound}.
\end{proof}

\begin{lemma}[Encoding projected vectors]\label{lem:encoding}
With the preceding notation,
\begin{equation}\label{eq:encoding}
N(U)\le |E(Q_U)\setminus E_0|+r+P+z
\le 4r+t_2+3h+P+z-6.
\end{equation}
\end{lemma}
\begin{proof}
We separate three exceptional families of vectors, always identifying opposite signs:
\begin{enumerate}
\item unit vectors $e_u$, $u\in U$: at most $r$ types;
\item pair indicators $e_u+e_v$ for pairs counted by $P$: at most $P$ types;
\item indicators $\ind_{A_F}$ for $|A_F|\ge3$ when $F$ has an edge-adjacent empty-type face: at most $z$ types.
\end{enumerate}
Every nonzero comparison involving an empty type belongs to one of these families: the other type has size one, two, or at least three, respectively.

For comparisons between nonempty types, equal labels give zero. Otherwise the dual edge representing the comparison survives as an edge of $Q_U$ joining the two type labels. If that simple edge belongs to $E_0$, its labels are $A=\{u,v\}$ or $A=\{u,v,w\}$ and the singleton $\{u\}\subset A$. The difference is respectively
\[
\ind_{\{u,v\}}-e_u=e_v,
\qquad\ind_{\{u,v,w\}}-e_u=e_v+e_w,
\]
which is already in the first or second exceptional family. In the second case $\{v,w\}$ is counted by $P$.

Every remaining nonzero comparison between nonempty types is represented by an edge of $E(Q_U)\setminus E_0$. This covers different singleton types, a singleton versus a type not containing it, a singleton versus a type of size at least four, and any two distinct nonsingleton types. The ordered type labels determine their difference up to sign, so one simple edge accounts for at most one remaining projected vector, regardless of original multiplicity. This proves the first inequality in \eqref{eq:encoding}. Substituting $R=r+t_2+t_3+h$ into \eqref{eq:quotientbound} gives the second.
\end{proof}

\begin{lemma}[Projected facial-difference types]\label{lem:projected-types}
For every nonempty $U\subseteq V(G)$,
\[
N(U)\le11|U|.
\]
More precisely, if $r=|U|\ge3$, then $N(U)\le11r-20$.
\end{lemma}
\begin{proof}
Suppose $r\ge3$. By Lemmas~\ref{lem:skeleton}--\ref{lem:encoding}, using $t_2,P\le e$, $z\le t_3+h$, and $e+h\le3r-6$, we obtain
\begin{equation}\label{eq:count1}
N(U)\le4r+2e+t_3+4h-6\le10r+t_3+2h-18.
\end{equation}
Alternatively, $t_2+z\le2e-3t_3-4h$ and $P\le e$ give
\begin{equation}\label{eq:count2}
N(U)\le4r+3e-3t_3-h-6\le13r-3t_3-4h-24.
\end{equation}
Taking twice \eqref{eq:count1} plus \eqref{eq:count2}, divided by three, yields
\[
N(U)\le11r-\frac{t_3}{3}-20\le11r-20.
\]
For $r=1$ or $2$, there are at most one or four nonzero vectors in $\{-1,0,1\}^r$ up to sign, respectively. Both are at most $11r$.
\end{proof}

\subsection{Independent groups and successive digits}
The projection estimate gives a density bound in a linear matroid. The next lemma converts a partition into independent groups into one modular labeling by handling one group at each digit.

\begin{lemma}[Successive modular digits]\label{lem:digits}
Let $L_1,\ldots,L_m$ be integer linear forms. Suppose their coefficient vectors modulo a prime $p$ can be partitioned into $t$ linearly independent sets. Then there is $b\in\{0,\ldots,p^t-1\}^n$ such that
\[
L_i(b)\not\equiv0\pmod{p^t}
\qquad\text{for all }i.
\]
\end{lemma}
\begin{proof}
Index the independent groups by $j=0,\ldots,t-1$ and construct
\[
b=\sum_{j=0}^{t-1}p^jz_j,
\qquad z_j\in\{0,\ldots,p-1\}^n.
\]
Suppose the partial sum $b_j=\sum_{h<j}p^hz_h$ is fixed. For each form in group $j$ still divisible by $p^j$, require
\[
L_i(z_j)\equiv1-\frac{L_i(b_j)}{p^j}\pmod p.
\]
The coefficient rows in group $j$ are linearly independent over $\FF_p$, so this linear system has a solution $z_j$. Every such form becomes nonzero modulo $p^{j+1}$. Forms already nonzero modulo a smaller power remain nonzero when higher digits are added. Induction proves the claim after $t$ digits.
\end{proof}

\begin{theorem}[Explicit universal upper bound]\label{thm:universal}
Every 3-connected plane graph $G$ satisfies
\[
\lambda(G)\le4\cdot3^{11}-1=708587.
\]
Consequently, $K_3(n)\le\min\{3n-6,708587\}$ and $8\le K_3^*<\infty$.
\end{theorem}
\begin{proof}
Reduce the facial-difference vectors \eqref{eq:rho} modulo three, retaining one representative from each nonzero projective class. Since all coefficients belong to $\{-1,0,1\}$, two such vectors are proportional over $\FF_3$ exactly when the corresponding integer vectors agree up to sign.

Let $S$ be any subset of representatives of rank $r$. Choose $r$ coordinate vertices $U$ such that projection of the row space spanned by $S$ onto $\FF_3^U$ is injective. All projected representatives remain nonzero and pairwise projectively distinct. Lemma~\ref{lem:projected-types} therefore gives
\[
|S|\le N(U)\le11r=11\rank_{\FF_3}(S).
\]
Edmonds' matroid partition theorem~\cite{edmonds} now partitions the representatives into at most eleven linearly independent sets. Applying Lemma~\ref{lem:digits} with $p=3$ yields a map $b:V(G)\to\{0,\ldots,3^{11}-1\}$ for which every represented facial difference is nonzero modulo $3^{11}$. The omitted signed duplicates satisfy the same condition.

The Four Color Theorem~\cite{fourcolor} supplies a proper four-coloring of $G$. Combining it with $b$ via Lemma~\ref{lem:modular-combination} gives the displayed constant. Proposition~\ref{prop:edge-count} supplies the additional bound in $n$, and \eqref{eq:tri-bounds} gives the lower bound on $K_3^*$.
\end{proof}

\begin{remark}[Quantitative role of the projected-type estimate]\label{rem:constant}
More generally, if some integer $t$ satisfies $N(U)\le t|U|$ for every nonempty $U$, the same argument yields $\lambda(G)\le4\cdot3^t-1$. Reducing the eleven independent groups to ten would give $236195$, and nine would give $78731$. These are conditional consequences, not improved bounds proved here. The relevant target is the rank-density of the projected facial-difference matroid. The value $708587$ is an explicit existence bound, not a claimed optimum or a practical search threshold.
\end{remark}

\begin{remark}\label{rem:prime-three}
The prime three distinguishes the signs of the original coefficient vectors; reducing modulo two would identify $1$ and $-1$. The final Four Color Theorem input is separate from projection, partition, and lifting. Replacing it by the elementary Five Color Theorem gives the alternative universal bound $5\cdot3^{11}-1=885734$.
\end{remark}

\section{Summary, open questions, and verification scope}\label{sec:outlook}
Table~\ref{tab:results} summarizes the numerical consequences, with $n=|V(G)|$ and $m=|E(G)|$.
\begin{table}[ht]
\centering\small
\caption{Summary of bounds. The final row records the classical diagonal-coloring bounds; the other bounds are proved in this paper.}\label{tab:results}
\begin{tabularx}{\textwidth}{@{}X p{4.3cm}@{}}
\toprule
Class or hypothesis & Bound\\
\midrule
All 3-connected plane graphs & $\lambda\le\min\{3n-6,708587\}$\\
All 3-connected plane graphs, arbitrary lists & $\laml\le m$\\
Cubic 3-connected plane graphs & $\lambda\le14$\\
Cubic graphs: one-class case, or two-class case satisfying \eqref{eq:modfour-condition} & $\lambda\le11$\\
Bipartite cubic 3-connected plane graphs & $\lambda\le9$\\
Fullerenes, arbitrary lists & $\laml\le7$\\
Triangulations with edge-disjoint all-high facial triangles & $\lambda\le8$ (sharp)\\
Triangulations with no all-high facial triangle & $\lambda\le7$\\
Triangulations with independent high vertices & $\chd\le7$, $\lambda\le6$ (sharp)\\
Unrestricted plane triangulations (classical) & $8\le K_{\mathrm{tri}}\le9$\\
\bottomrule
\end{tabularx}
\end{table}

The principal open questions are whether every plane triangulation admits a diagonal nine-coloring, whether the fullerene common-alphabet bound can be reduced to three, and how far the universal constant can be lowered. These directions involve different obstacles. The modular lifting method produces a common integer alphabet, not an arbitrary-list theorem. Conversely, an eight-element arbitrary-list pentagon extension does not by itself yield the conjectured four-label fullerene theorem. In the cubic proof, we use modulus five to handle the two-class cases not covered by the sufficient condition \eqref{eq:modfour-condition}; we do not claim that modulus four is impossible in every such case. In the universal argument, reducing the projected-type rank density is the main quantitative objective.

\subsection{Proof dependencies and supplementary verification}
The coefficient identity \eqref{eq:coefficient} has a symbolic, self-contained proof in Appendix~\ref{app:coefficient}. Optional verification scripts are supplied in the \texttt{supplementary/} directory accompanying this revision. They use only the Python standard library. The script \texttt{verify\_symbolic.py} checks the polynomial identities, alternating signs, and tableau counts; \texttt{verify\_coefficient\_direct.py} independently extracts the target coefficient modulo seven from the original factors. The script \texttt{verify\_balancing.py} checks the finite branches of Lemma~\ref{lem:balance-five}. Numerical graph examples are checked by \texttt{verify\_examples.py}, which also runs \texttt{verify\_projections.py} on specified embedded graphs and selected vertex subsets.

The accompanying README identifies this manuscript version, gives the commands and expected outputs, and describes the coverage and limitations of each check. The recorded run report gives the interpreter, elapsed times, and exit codes; a SHA-256 manifest identifies the distributed files. These finite computations provide reproducible checks of the displayed algebra and examples. They do not prove the local constraint classification or the projected-type bounds for all graphs, and they are not premises of the stated theorems.

Some imported theorem-level ingredients have their own computational history. The Four Color Theorem~\cite{fourcolor} is used in Theorems~\ref{thm:transversal}, \ref{thm:cubic}, and~\ref{thm:universal}. Other named inputs are Borodin's $1$-planar six-color theorem, list Brooks coloring, the Combinatorial Nullstellensatz, cyclic edge-connectivity of fullerenes, and matroid partition. Optional computational checks do not constitute additional hypotheses in any theorem.

\appendix
\section{A symbolic proof of the pentagon coefficient}\label{app:coefficient}
All computations in this appendix take place in $k=\FF_7$, except where an integer count is explicitly indicated. Let
\[
\Delta(x)=\prod_{0\le i<j\le4}(x_i-x_j),\qquad
S=\sum_{i=0}^4x_i,\qquad y_i=x_i+x_{i+1},
\]
with indices modulo five. Removing $\prod_i x_i$ from \eqref{eq:Qtilde}, the coefficient to be computed is
\begin{equation}\label{eq:A1}
A=[x_0^6x_1^6x_2^6x_3^6x_4^6]\,\Delta(x)
\prod_{i=0}^4 y_i^3\prod_{i=0}^4(S-y_i).
\end{equation}
We prove $A=2$ through symbolic identities and elementary counting. In particular, the evaluations used below determine a previously established two-dimensional space of possible symmetric quotients; they are not substituted for a polynomial identity.

\subsection{The top coefficient under linear substitutions}
\begin{lemma}\label{lem:top-invariance}
Let $P\in\FF_p[X_1,\ldots,X_d]$ be homogeneous of degree $d(p-1)$, and let $M\in GL_d(\FF_p)$. Then
\[
[X_1^{p-1}\cdots X_d^{p-1}]P(MX)
=(\det M)^{p-1}[X_1^{p-1}\cdots X_d^{p-1}]P(X).
\]
In particular, for $p=7$ the coefficient is invariant under invertible linear substitutions.
\end{lemma}
\begin{proof}
Work in $\FF_p[X_1,\ldots,X_d]/(X_1^p,\ldots,X_d^p)$. Linear substitutions preserve this ideal, since in characteristic $p$,
\[
\Bigl(\sum_i a_iX_i\Bigr)^p=\sum_i a_i^pX_i^p.
\]
Its component of degree $d(p-1)$ is one-dimensional, generated by $\prod_iX_i^{p-1}$. Diagonal scaling acts by $(\det M)^{p-1}$, permutation matrices act by $(\operatorname{sgn}\pi)^{p-1}$, and transvections act trivially because any other term has an exponent at least $p$. These elementary matrices generate $GL_d$, proving the formula. Fermat's identity yields $(\det M)^6=1$ when $p=7$.
\end{proof}

Put $w=S/5$ and $z_i=x_i-w$, so $\sum_i z_i=0$. Take independent coordinates $w,z_0,z_1,z_2,z_3$, with $z_4=-z_0-z_1-z_2-z_3$. This is an invertible linear change of variables, and $\Delta(x)=\Delta(z)$. Set
\[
t_i=z_i+z_{i+1},\qquad a_j=e_j(t_0,\ldots,t_4),
\]
where $e_j$ denotes the $j$th elementary symmetric polynomial. Since $a_1=0$, the products
\[
C(w)=\prod_i(2w+t_i),\qquad D(w)=\prod_i(3w-t_i)
\]
have the sparse expansions
\begin{align*}
C(w)&=32w^5+8a_2w^3+4a_3w^2+2a_4w+a_5,\\
D(w)&=243w^5+27a_2w^3-9a_3w^2+3a_4w-a_5.
\end{align*}
Collecting the coefficient of $w^6$ yields
\begin{equation}\label{eq:A2}
\begin{split}
R(z):=[w^6]C(w)^3D(w)
={}&4a_4a_5^2+a_2^2a_5^2+a_2a_3a_4a_5\\
&{}+a_3^3a_5+3a_3^2a_4^2.
\end{split}
\end{equation}
For transparency, the six possible $w$-degree multisets in the four factors $C,C,C,D$, together with their total integer contributions before reduction modulo seven, are shown in Table~\ref{tab:coeff-contrib}.
\begin{table}[ht]
\centering\small
\caption{Contributions to $[w^6]C(w)^3D(w)$; reductions use $k=\FF_7$.}\label{tab:coeff-contrib}
\begin{tabular}{@{}cccc@{}}
\toprule
$w$-degrees & Factor in $a_j$ & Integer contribution & Modulo $7$\\
\midrule
$(5,1,0,0)$ & $a_4a_5^2$ & $1362$ & $4$\\
$(3,3,0,0)$ & $a_2^2a_5^2$ & $456$ & $1$\\
$(3,2,1,0)$ & $a_2a_3a_4a_5$ & $624$ & $1$\\
$(3,1,1,1)$ & $a_2a_4^3$ & $504$ & $0$\\
$(2,2,2,0)$ & $a_3^3a_5$ & $-496$ & $1$\\
$(2,2,1,1)$ & $a_3^2a_4^2$ & $-144$ & $3$\\
\bottomrule
\end{tabular}
\end{table}
For example, the first integer contribution is $3\cdot243\cdot2+3\cdot3\cdot32-6\cdot32\cdot2=1362$, according as the $D$ factor has $w$-degree five, one, or zero. The other rows follow by the same multiplication of the two displayed sparse polynomials, and the six degree patterns exhaust all ways for four selected degrees to sum to six. By Lemma~\ref{lem:top-invariance},
\begin{equation}\label{eq:A3}
A=[z_0^6z_1^6z_2^6z_3^6]\,\Delta(z)R(z).
\end{equation}

\subsection{An alternating identity on the centered hyperplane}
The alternating quotient below belongs to the two-dimensional space spanned by $e_2(z)^2$ and $e_4(z)$ on the hyperplane $\sum_i z_i=0$. A specialization on which these basis elements remain linearly independent therefore determines the two coefficients.

For a polynomial in the centered coordinates, define
\[
\Alt R(z)=\sum_{\sigma\in S_5}\operatorname{sgn}(\sigma)
R(z_{\sigma(0)},\ldots,z_{\sigma(4)}).
\]
\begin{lemma}\label{lem:alternant-identity}
The polynomial $R$ in \eqref{eq:A2} satisfies
\begin{equation}\label{eq:A4}
\Alt R(z)=3\Delta(z)\bigl(e_2(z)^2+e_4(z)\bigr)
\qquad\text{when }\sum_{i=0}^4z_i=0.
\end{equation}
\end{lemma}
\begin{proof}
Extend centered polynomials to five free variables by replacing $z_i$ with $z_i-(\sum_jz_j)/5$. Any alternating polynomial is divisible by $\Delta$, and the quotient here is symmetric, centered, and of degree four. The elementary symmetric-polynomial theorem therefore gives
\[
\Alt R(z)=\Delta(z)\bigl(\alpha e_2(z)^2+\beta e_4(z)\bigr)
\]
for some $\alpha,\beta\in k$. It suffices to specialize
\[
z=(-t,-1,0,1,t),
\]
because $e_2(z)^2=(t^2+1)^2$ and $e_4(z)=t^2$ are linearly independent polynomials in $t$.

The summand $R$ is invariant under rotation and reversal of its cyclic order. This dihedral group has ten elements, all even permutations of five positions. Thus we may sum over unoriented cyclic orders and multiply by ten.

If a cyclic order includes both edges joining the opposite pairs, two pair sums vanish, giving $a_4=a_5=0$ and therefore $R=0$. If it contains exactly one such edge, the pair sums are $\{0,b,-b,c,-c\}$ for suitable $b,c$. Then $a_3=a_5=0$ and again $R=0$. Consequently only orders avoiding both opposite-pair edges contribute. Deleting the zero vertex from such an order leaves a Hamilton path of the four-cycle obtained from $K_4$ by deleting the opposite-pair matching. Each unoriented Hamilton path is obtained by deleting one of the four cycle edges, so there are exactly four.

Choose the representative cyclic order $(-t,-1,0,t,1)$, an odd permutation of the displayed base tuple. Its consecutive pair sums are
\[
(-t-1,-1,t,t+1,1-t).
\]
Their elementary symmetric polynomials are
\begin{align*}
a_2&=-2t^2-t-2,& a_3&=t(t-1),\\
a_4&=(t+1)^2(t^2-t+1),& a_5&=-t(t-1)(t+1)^2.
\end{align*}
Substituting these into \eqref{eq:A2} and multiplying in $k[t]$ gives
\[
R_0(t)=t^2(t^2-1)^2(-t^6-t^5+4t^3-t-1).
\]
The four contributing cyclic orders are listed in Table~\ref{tab:cyclic-signs}. Their permutation signs are measured relative to the base tuple $(-t,-1,0,1,t)$. Swapping the pair $-t,t$ replaces $t$ by $-t$ in $R_0$ and reverses the sign. Swapping both opposite pairs negates all five entries and preserves the sign; since $R$ has even degree, it leaves the value of $R$ unchanged. These observations determine all four entries of the table.
\begin{center}
\begin{minipage}{\textwidth}
\centering\small
\captionsetup{hypcap=false}
\captionof{table}{The four contributing unoriented cyclic orders in the specialization $z=(-t,-1,0,1,t)$.}\label{tab:cyclic-signs}
\begin{tabular}{@{}ccc@{}}
\toprule
Representative cyclic order & Permutation sign & Signed contribution\\
\midrule
$(-t,-1,0,t,1)$ & $-1$ & $-R_0(t)$\\
$(t,-1,0,-t,1)$ & $+1$ & $+R_0(-t)$\\
$(-t,1,0,t,-1)$ & $+1$ & $+R_0(-t)$\\
$(t,1,0,-t,-1)$ & $-1$ & $-R_0(t)$\\
\bottomrule
\end{tabular}
\end{minipage}
\end{center}
Their signed sum is therefore $-2(R_0(t)-R_0(-t))$. Including the factor ten gives
\begin{align*}
\Alt R(-t,-1,0,1,t)
&=-20\bigl(R_0(t)-R_0(-t)\bigr)\\
&=5t^3(t^2-1)^2(t^4+3t^2+1)\qquad\text{in }k[t].
\end{align*}
On the other hand,
\[
\Delta(-t,-1,0,1,t)=4t^3(t^2-1)^2,
\qquad e_2(z)^2+e_4(z)=t^4+3t^2+1.
\]
Since $5/4=3$ in $k$, it follows that $\alpha=\beta=3$, proving \eqref{eq:A4}.
\end{proof}

Permutations of the five centered coordinates act invertibly on the four independent $z$-coordinates. The coefficient functional in \eqref{eq:A3} is invariant under this action by Lemma~\ref{lem:top-invariance}. As $5!=120=1$ in $k$, replacing $R$ by its alternation in \eqref{eq:A3} yields
\begin{equation}\label{eq:A5}
A=3B,\qquad
B=[z_0^6z_1^6z_2^6z_3^6]\,\Delta(z)^2\bigl(e_2(z)^2+e_4(z)\bigr).
\end{equation}

To return to five free variables, apply Lemma~\ref{lem:top-invariance} to
$S^6\Delta(x)^2\bigl(e_2(z)^2+e_4(z)\bigr)$, with $z_i=x_i-S/5$. Since $S=5w$ and $5^6=1$ in $k$, its top coefficient equals $B$. The elementary shift identities give
\[
e_2(z)=e_2(x)+S^2,\qquad
e_4(z)=e_4(x)+Se_3(x)-S^2e_2(x)+3S^4.
\]
In $k[x_0,\ldots,x_4]/(x_0^7,\ldots,x_4^7)$, we have $S^7=\sum_i x_i^7=0$. After multiplying by $S^6$, all shift terms vanish. Therefore
\begin{equation}\label{eq:A6}
B=[x_0^6x_1^6x_2^6x_3^6x_4^6]\,
\Delta(x)^2 S^6\bigl(e_2(x)^2+e_4(x)\bigr).
\end{equation}

\subsection{The final coefficient from elementary alternants}
For a partition $\mu$ with at most five parts, define the Schur polynomial by the alternant quotient
\[
s_\mu(x)=\frac{\det\bigl(x_i^{\mu_j+5-j}\bigr)_{0\le i\le4,\,1\le j\le5}}{\Delta(x)}.
\]
These symmetric polynomials have integer coefficients; see Macdonald~\cite{macdonald} or Stanley~\cite{stanley} for background. We use two elementary alternant facts, each briefly justified below.

First, for any homogeneous symmetric polynomial $f$ of degree ten,
\begin{equation}\label{eq:A7}
[x_0^6\cdots x_4^6]\,\Delta(x)^2f(x)
=5!\,[s_{(2,2,2,2,2)}]f,
\end{equation}
where the bracket on the right denotes the coefficient in the Schur expansion. Indeed, if $f=s_\mu$, the numerator alternant has strictly decreasing exponents. Combining it with a permutation of exponents $0,1,2,3,4$ from the second Vandermonde to obtain exponent six on every variable forces its exponent list to be $(6,5,4,3,2)$, hence $\mu=(2,2,2,2,2)$. For that partition all $5!$ complementary permutation pairs have positive sign, because reversal of five entries is even.

Second, multiplication by $e_1=S$ adds one cell to a partition:
\begin{equation}\label{eq:A8}
e_1s_\mu=\sum_{\nu:\,\nu/\mu\text{ is one cell}}s_\nu.
\end{equation}
This follows by multiplying the numerator alternant by $\sum_i x_i$. The resulting terms raise one exponent by one; repeated exponents make an alternant zero, and the surviving sequences are exactly the one-cell extensions.

The two-column determinant identities also give
\[
s_{(2,2)}=e_2^2-e_1e_3,\quad
s_{(2,1,1)}=e_1e_3-e_4,\quad
s_{(1,1,1,1)}=e_4.
\]
Consequently,
\begin{equation}\label{eq:A9}
e_2^2+e_4=s_{(2,2)}+s_{(2,1,1)}+2s_{(1,1,1,1)}.
\end{equation}

Iterating \eqref{eq:A8} six times counts standard Young tableaux of the skew shape obtained by removing each starting partition from the rectangle $(2^5)$. Rotate the complement by $180^\circ$ and reverse the entries to identify these counts with those of straight shapes. Transposing these shapes preserves the number of standard tableaux. Table~\ref{tab:young} records the resulting counts and the multiplicities from \eqref{eq:A9}.
\begin{table}[ht]
\centering\small
\caption{Young shapes contributing to the rectangular Schur coefficient.}\label{tab:young}
\begin{tabular}{@{}cccc@{}}
\toprule
Starting $\mu$ & Rotated complement & Transpose & Multiplicity\\
\midrule
$(2,2)$ & $(2,2,2)$ & $(3,3)$ & 1\\
$(2,1,1)$ & $(2,2,1,1)$ & $(4,2)$ & 1\\
$(1,1,1,1)$ & $(2,1,1,1,1)$ & $(5,1)$ & 2\\
\bottomrule
\end{tabular}
\end{table}

Each skew shape has six cells. On rotation, replacing every entry $j$ with $7-j$ restores the increasing row and column conditions; transposition preserves the number of standard tableaux. For a two-row shape $(a,b)$, a standard ballot reflection argument gives
\[
f^{(a,b)}=\binom{a+b}{b}-\binom{a+b}{b-1}.
\]
Choose the entries of the second row; sequences violating the ballot condition are counted by reflecting at their first prefix with more lower-row than upper-row entries. The relevant numbers are
\[
\binom63-\binom62=5,\qquad
\binom62-\binom61=9,\qquad
\binom61-\binom60=5.
\]
Combining \eqref{eq:A7}--\eqref{eq:A9} therefore yields
\[
B=5!\,(5+9+2\cdot5)=120\cdot24=3\qquad\text{in }\FF_7.
\]
Finally, \eqref{eq:A5} gives $A=3\cdot3=2$ in $\FF_7$. By \eqref{eq:A1}, this is exactly the nonzero coefficient asserted in \eqref{eq:coefficient}.

\section{A separate elementary fullerene bound}\label{app:elementary}
This appendix gives a thirteen-element list theorem without using the Combinatorial Nullstellensatz or the coefficient identity. It is weaker than Theorem~\ref{thm:fullerene-list}, but provides a direct alternative.

\begin{lemma}\label{lem:B1}
If a vertex $x$ has an unlabeled neighbor $y$, and the face at $x$ not containing $xy$ has length $c$, then at most $c+2\le8$ constraints can become determined when $x$ is labeled in a fullerene.
\end{lemma}
\begin{proof}
If the two faces at $xy$ have lengths $a,b$, at most $a+b+c-3$ factor occurrences involve $x$. The edge factor at $xy$, as well as the $a-3$ and $b-3$ facial factors indexed by edges incident with neither endpoint, still involve the unlabeled $y$. Subtracting these $a+b-5$ deferred occurrences leaves at most $c+2$.
\end{proof}

\begin{lemma}\label{lem:B2}
For finite sets $A,B\subset\R$ of sizes $s,t\ge2$, the restricted sumset $\{a+b:a\in A,\ b\in B,\ a\ne b\}$ has at least $s+t-3$ elements.
\end{lemma}
\begin{proof}
Write $A=\{a_1<\cdots<a_s\}$ and $B=\{b_1<\cdots<b_t\}$. The chain
\[
a_1+b_1<\cdots<a_1+b_t<a_2+b_t<\cdots<a_s+b_t
\]
contains $s+t-1$ distinct sums. At most one entry in each of its two monotone segments has equal summands. Removing those at most two entries leaves at least $s+t-3$ admissible sums.
\end{proof}

\begin{theorem}\label{thm:B3}
Every fullerene satisfies $\laml(G)\le12$. If it has two adjacent pentagons, then $\laml(G)\le11$.
\end{theorem}
\begin{proof}
Choose an edge $uv$ of a pentagon and leave its endpoints unlabeled. Process all remaining vertices in nonincreasing distance from $\{u,v\}$. Every processed vertex has an uncolored neighbor, so Lemma~\ref{lem:B1} permits a choice from a thirteen-element list.

Let $P,Q$ be the faces at $uv$, and $R,T$ the third faces at $u,v$, respectively. Write
\[
|P|=a=5,\qquad |Q|=b\le6,\qquad |R|=c\le6,\qquad |T|=d\le6.
\]
The faces $R,T$ are not adjacent: otherwise their four corresponding dual vertices would induce a $K_4$, which in a dual triangulation with more than four vertices forces a separating triangle, contrary to Lemma~\ref{lem:pentagon-neighbors}.

Write $x=c(u)$ and $y=c(v)$. Of the remaining constraints, at most $c+2$ restrict $x$ alone and at most $d+2$ restrict $y$ alone. There is also the edge restriction $x\ne y$ and at most
\[
(a-3)+(b-3)\le5
\]
forbidden values of $x+y$. This list is exhaustive: only $P,Q$ contain both endpoints; comparisons with $R,T$ are one-variable restrictions; and $R,T$ are not adjacent.

After deleting the one-variable forbidden values, each thirteen-element list retains at least five values. Lemma~\ref{lem:B2} produces at least $5+5-3=7$ admissible restricted sums, more than the at most five forbidden values of $x+y$. This completes the labeling. If $a=b=5$, twelve-element lists still retain four values per endpoint, yielding at least five restricted sums against at most four forbidden sums. The same greedy stage applies, proving $\laml(G)\le11$ in this case.
\end{proof}
\section*{Acknowledgements}
This project is  supported by the National Natural Science Foundation of China (No. 12271246). 
\section*{AI Declaration}

We used ChatGPT to assist with the exploration of proof strategies and the preparation of this manuscript.  All mathematical arguments and proofs in the final manuscript were written and verified by the author.

\end{document}